\documentclass[10pt]{amsart}

\usepackage{amssymb,eucal,amsmath,amsthm,graphicx,mathrsfs,times}
\usepackage[latin1]{inputenc}

\newcommand{\mt}[1]{\text{\rm #1}}    % text in math mode
\newcommand{\comment}[1]{}            % usage: {*** ignored text ***}
\newcommand{\Proof}{\begin{proof}[\textsc{Beweis}] \;}
\newcommand{\set}[1]{\{#1\}}
\newcommand{\suchthat}{\;|\;}         % for sets: {x | x>z}
\newcommand{\restrict}{\,|}           % restriction of maps: f|M
\newcommand{\compose}{\circ}          % composition of maps: f°g
\newcommand{\define}{\;{\rm :=}\;}
\newcommand{\invdef}{\;{\rm =:}\;}

\newcommand{\without}{\mathord{\setminus}}

\newcommand{\pr}{\mt{pr}}             % projection 2
\newcommand{\R}{\mathbb{R}}           % real numbers
\newcommand{\N}{\mathbb{N}}           % natural numbers
\newcommand{\Z}{\mathbb{Z}}           % integers
\newcommand{\leer}{\varnothing}       % empty set
\newcommand{\mfbd}{\partial}          % manifold boundary
\newcommand{\eval}[2]{\langle #1,#2 \rangle}     % scalar product
\newcommand{\abs}[1]{\lvert#1\rvert}  % absolute value
\newcommand{\norm}[1]{\lVert#1\rVert} % norm
\newcommand{\diag}{\mt{diag}}         % diagonal matrix
\newcommand{\Gr}{\mt{Gr}}             % Grassmannian
\newcommand{\eps}{\varepsilon}        % standard epsilon
\newcommand{\ot}{\text{$\gets$}}      % C^r(M\ot E): C^r sections in a bundle E \to M
\newcommand{\ocinterval}[2]{\mathord{\mathopen{]}#1,#2\mathclose{]}}}  % left-open interval
\swapnumbers
\theoremstyle{definition}
\newtheorem{remark}{Remark}[section]

\newtheorem{remarks}[remark]{Remarks}

\newtheorem{generalconventions}[remark]{General conventions}

\newtheorem{definition}[remark]{Definition}

\newtheorem{notation}[remark]{Notation}

\newtheorem{distributionspace}[remark]{The space of distributions}
\newtheorem{conncomp}[remark]{Connected components}

\theoremstyle{plain}
\newtheorem{theorem}[remark]{Theorem}

\newtheorem{lemma}[remark]{Lemma}

\newtheorem{keylemma}[remark]{Key lemma}
\newtheorem{corollary}[remark]{Corollary}
\newtheorem{proposition}[remark]{Proposition}

\DeclareMathOperator{\im}{\mt{im}}

\newcommand{\Sym}{\mt{Sym}}       % symmetric bilinear forms
\DeclareMathOperator{\trace}{tr}
\DeclareMathOperator{\scal}{scal}
\DeclareMathOperator{\Ric}{Ric}
\DeclareMathOperator{\Riem}{Riem}

\DeclareMathOperator{\Hess}{Hess}
\DeclareMathOperator{\divergence}{div}
\DeclareMathOperator{\laplace}{\Delta}
\DeclareMathOperator{\spann}{span}

\DeclareMathOperator{\Lin}{Lin}

\renewcommand{\labelenumi}{(\roman{enumi})}
\renewcommand{\bar}{\overline}

\newcommand{\Quiem}{\mt{Quiem}}              % the Quiemann tensor
\DeclareMathAlphabet{\mathpzc}{OT1}{pzc}{m}{it}

\newcommand{\stre}{\text{\tt stretch}}

\newcommand{\switch}{\text{\tt switch}}

\DeclareMathOperator{\Distr}{Distr}
\newcommand{\Twist}{\mathpzc{Tw}}
\newcommand{\Symst}{\mathpzc{Sw}}
\DeclareMathOperator{\Metr}{Metr}

\newcommand{\LeviCivita}{Levi-Civita}

\newcommand{\role}{role}

\newcommand{\z}{s}                     % realise/realize etc.
\newcommand{\ou}{ou}                   % neighbourhood/neighborhood etc.
\newcommand{\re}{re}                   % centre/center etc.
\begin{document}

%\selectlanguage{english}

\title[Nonexistence of space foliations and the dominant energy condition]{Nonexistence of spacelike foliations and\\
the dominant energy condition in Lorentzian geometry}
\author{Marc Nardmann}
\address{Department of Mathematics, University of Regensburg}
\email{Marc.Nardmann\@@mathematik.uni-regensburg.de}
\thanks{This work was supported by the priority programme \emph{Globale Differentialgeometrie} of the Deutsche For\-schungs\-gemeinschaft (DFG)}

\begin{abstract}
We show that many Lorentzian manifolds of dimension $\geq3$ do not admit a spacelike codimension-one foliation, and that almost every manifold of dimension $\geq3$ which admits a Lorentzian metric at all admits one which satisfies the dominant energy condition and the timelike convergence condition. These two seemingly unrelated statements have in fact the same origin.

\smallskip\noindent
We also discuss the problem of topology change in General Relativity. A theorem of Tipler says that topology change is impossible via a spacetime cobordism whose Ricci curvature satisfies the strict lightlike convergence condition. In his theorem, the boundary of the cobordism is required to be spacelike. We show that topology change with the strict lightlike convergence condition and also the dominant energy condition is possible in many cases when one requires instead only that there exists a timelike vector field which is transverse to the boundary.
\end{abstract}

\maketitle

\setcounter{section}{-1}

\section{Introduction}

We consider two questions in Lorentzian geometry which seem unrelated at first sight:
\begin{itemize}
\item
Does every Lorentzian manifold admit a spacelike codimension-one foliation?
\item
When a manifold admits a Lorentzian metric at all, does it admit one which satisfies the dominant energy condition?
\end{itemize}

The first question has actually more to do with differential topology than Lorentzian geometry. In particular, whether a Lorentzian manifold admits a spacelike codimension-one foliation depends only on the conformal class of the metric. As we will see, many Lorentzian manifolds of dimension $\geq3$ do not admit a spacelike codimension-one foliation.

\smallskip
In contrast, the second question is a Ricci curvature problem. The dominant energy condition, which plays an important {\role} in General Relativity, is a Ricci nonnegativity condition which depends on a ``cosmological'' constant $\Lambda$ (the definition is reviewed in Section \ref{backgroundE}). It turns out that the answer to the second question is always \emph{yes} in dimension $\geq5$, and it is \emph{yes} in dimensions $3$ and $4$ under mild assumptions.

\smallskip
Somewhat surprisingly, there is a close connection between the two questions. Let us discuss nonexistence of spacelike codimension-one foliations first.

\smallskip
Most research in Lorentzian Geometry and General Relativity deals with metrics which have nice causality properties like global hyperbolicity or stable causality. A Lorentzian manifold with these properties admits a smooth function with timelike gradient (\cite{BernalSanchez2005}, Theorem 1.2) and thus, in particular, a smooth foliation by spacelike codimension-one submanifolds (namely by the level sets of the time function).

\smallskip
The question arises what happens when we drop the causality condition: Does \emph{every} Lorentzian manifold admit a smooth spacelike codimension-one foliation? (Note that we do not demand that the leaves be submanifolds; they might be dense, for instance, as in the case of the foliation of the torus $\R^2/\Z^2$ by lines of irrational slope.) Every Lorentzian manifold admits a smooth spacelike corank-one sub vector bundle of the tangent bundle, but it might happen that no such subbundle exists which is \emph{integrable}, i.e., which is the tangent bundle to some foliation on $M$.

\smallskip
Clearly every point in a Lorentzian manifold has a neighb{\ou}rhood which admits a spacelike codimension-one foliation. Every two-dimensional Lorentzian manifold admits a spacelike (co)dimension-one foliation simply because every rank-one subbundle of the tangent bundle is integrable. In higher dimensions, nonintegrable corank-one subbundles always exist; so the question whether one of the spacelike corank-one subbundles is integrable becomes nontrivial. By a famous theorem of W.\ Thurston (\cite{Thurston1976}, Theorem 1), there are no purely topological obstructions to the existence of codimension-one foliations: Every connected component of the space of corank-one distributions on a given manifold contains an integrable distribution. (The word \emph{distribution} is always used in the sense of differential topology here: a \emph{$k$-plane distribution}, or synonymously: \emph{rank-$k$ distribution}, on a manifold is a rank-$k$ sub vector bundle of the tangent bundle.) This implies that every connected component of the space of Lorentzian metrics on a given manifold contains a metric which admits a spacelike codimension-one foliation.

\smallskip
In spite of this, spacelike codimension-one foliations do not exist for every metric. Although we are mainly interested in the Lorentzian case $q=1$, we state this result for pseudo-Riemannian metrics of arbitrary index $q$. Let us call a spacelike codimension-$q$ foliation on a pseudo-Riemannian manifold of index $q$ a \emph{space foliation}, for simplicity.

\begin{theorem} \label{main0}
Let $(M,g)$ be an $n$-dimensional pseudo-Riemannian manifold of index $q\in\set{1,\dots,n-2}$ (e.g.\ a Lorentzian manifold of dimension $n\geq3$). Let $A\neq M$ be a closed subset of $M$. Then there exists a metric $g'$ of index $q$ on $M$ such that
\begin{itemize}
\item
$g=g'$ on $A$;
\item
every $g$-timelike vector in $TM$ is $g'$-timelike;
\item
$M\without A$ does not admit any codimension-$q$ foliation none of whose tangent vectors is $g'$-timelike; in particular, $(M,g')$ does not admit any space foliation.
\end{itemize}
\end{theorem}

Note that the nonexistence of spacelike codimension-one foliations is not a matter of complicated manifold or bundle topology: We can e.g.\ take $(M,g)$ to be Minkowski space and $A$ to be the complement of an arbitrarily small open ball.

\smallskip
Theorem \ref{main0} is not particularly deep; it follows from elementary facts of differential topology, as we will see in Section \ref{main0proof}. Things become more complicated and interesting when we ask to which extent nonexistence of space foliations is related to curvature properties, e.g.\ to the dominant energy condition (which all ``physically reasonable'' metrics in General Relativity are assumed to satisfy). Are there metrics which do not admit a space foliation but satisfy the dominant energy condition?

\smallskip
Surprisingly, it turns out that not only the answer is \emph{yes}; but that both properties, nonexistence of space foliations and the dominant energy condition, have a tendency to hold simultaneously (at least within certain $1$-parameter families of metrics; cf.\ Section \ref{sectiondominant} for details). Intuitively speaking, if one deforms a Lorentzian metric in a natural way such that the ``energy dominance'' becomes stronger and stronger, then at some point space foliations cease to exist. Conversely, if one deforms a metric by squeezing the set of spacelike vectors in the tangent bundle in such a way that space foliations cease to exist, then the energy dominance has a tendency to become stronger, so that eventually the dominant energy condition holds.

\smallskip
A similar link exists when we replace the dominant energy condition by the timelike convergence condition. (Recall that a Lorentzian metric $g$ satisfies the \emph{timelike convergence condition} iff $\Ric_g(v,v)\geq0$ holds for all timelike vectors $v\in TM$.)

\smallskip
In Section \ref{specialcase}, we study this link in the simplest special case, a certain $1$-parameter family of Lorentzian metrics on $\R^n$. This special case has the advantage over our later more general considerations that one can also discuss the behavi{\ou}r of geodesics. In particular, the special case shows that metrics without space foliation can be geodesically complete. Moreover, the unavoidable failure of causality conditions for such metrics can be seen very explicitly here. The precise results are as follows.

\begin{definition} \label{maindef}
For $c\in\R_{>0}$, we consider the following frame $(e_0,e_1,e_2)$ of the vector bundle $T\R^3$, induced by the vector fields $\partial_i=\frac{\partial}{\partial x_i}$ coming from the standard coordinates $(x_0,x_1,x_2)$ on $\R^3$:
\begin{align*}
e_0 &\define \frac{1}{c}\partial_0 \;\;,  & e_1 &\define \partial_1-x_2\partial_0 \;\;,  & e_2 &\define \partial_2 \;\;.
\end{align*}
We define the Lorentzian metric $g^c_3$ on $\R^3$ by declaring $(e_0,e_1,e_2)$ to be an orthonormal frame for which the vector field $e_0$ is timelike. I.e., $g^c_3(e_i,e_j)=\eps_i\delta_{ij}$ for $i,j\in\set{0,1,2}$, where $\eps_0=-1$ and $\eps_i=1$ if $i>0$, and where $\delta_{ij}$ denotes the Kronecker symbol. For $n\geq3$, we define the Lorentzian metric $g^c_n$ on $\R^n = \R^3\times\R^{n-3}$ to be the product metric of $g^c_3$ with the euclidean (Riemannian) metric on $\R^{n-3}$.
\end{definition}

\begin{remark} \label{mainremark}
Let $n\geq3$. The following statements hold for all $c\in\R_{>0}$:
\begin{itemize}
\item
The diffeomorphism $\varphi_c\colon\R^n\to\R^n$ given by
\[
(x_0,x_1,x_2,\dots,x_{n-1})\mapsto (c^2x_0,cx_1,cx_2,\dots,cx_{n-1})
\]
satisfies $g^c_n = \frac{1}{c^2}\varphi_c^\ast(g^1_n)$; thus $c^2 g^c_n$ is isometric to $g^1_n$.
\item
$g^c_n$ is geodesically complete.
\item
$g^c_n$ has no closed causal geodesic. But for every $p\in\R^n$, there is a $1$-parameter family of closed spacelike geodesics through $p$.
\item
Let $I\subseteq\R$ be a compact interval. Then every $g^c_3$-timelike path $I\to\R^3$ admits an extension to a closed $g^c_3$-timelike path. (We consider $C^r$ timelike paths here, where $r$ is any element of $\N_{\geq1}\cup\set{\infty}$; the statement holds for each $r$.) Hence there exist for every $p\in\R^n$ closed $g^c_n$-timelike paths through $p$.
\item
$g^c_n$ satisfies the timelike convergence condition and, for all $\Lambda\leq\frac{c^2}{4}$, the dominant energy condition with respect to the cosmological constant $\Lambda$.
\end{itemize}
\end{remark}

\begin{theorem} \label{mainRn}
Let $n\geq3$. To every $\Lambda\in\R$ and every nonempty open set $U\subset\R^n$, there exists a number $c_0>0$ such that for all $c\geq c_0$, the metric $g^c_n$ has the following properties:
\begin{itemize}
\item
$U$ does not admit a codimension-one foliation none of whose tangent vectors is $g^c_n$-timelike; in particular, $U$ does not admit a $g^c_n$-space foliation.
\item
$g^c_n$ satisfies the dominant energy condition with cosmological constant $\Lambda$.
\end{itemize}
\end{theorem}

After this special case, we are now going to discuss the situation on general manifolds:

\begin{theorem} \label{maindom}
Let $(M,g)$ be a connected Lorentzian manifold of dimension $n\geq4$, let $K$ be a compact subset of $M$, let $\Lambda\in\R$. If $n=4$, assume that $(M,g)$ is time- and space-orientable, and that either $M$ is noncompact, or compact with intersection form signature divisible by $4$. Then there exists a Lorentzian metric $g'$ on $M$ such that
\begin{itemize}
\item
every $g$-causal vector in $TM$ is $g'$-timelike;
\item
$g'$ satisfies the timelike convergence condition on the set $K$;
\item
$g'$ satisfies the dominant energy condition with cosmological constant $\Lambda$ on $K$;
\item
$M$ does not admit any codimension-one foliation none of whose tangent vectors is $g'$-timelike; in particular, $(M,g')$ does not admit a space foliation.
\end{itemize}
\end{theorem}

This theorem generali{\z}es to dimension $3$ when one assumes that $M$ is orientable and admits a $g$-spacelike contact structure; cf.\ Theorem \ref{maindom3} below. When one assumes only that $M$ is orientable, then Theorem \ref{maindom} holds with the statement ``\emph{every $g$-causal vector in $TM$ is $g'$-timelike}'' replaced by the weaker statement ``\emph{$g'$ lies in the same connected component of the space of Lorentzian metrics as $g$}''.

\medskip
Closely related to Theorem \ref{maindom}, we get new insights into the classical problem of ``topology change'' in General Relativity. Let us say that a Lorentzian manifold $(M,g)$ satisfies the \emph{strict lightlike convergence condition} iff $\Ric_g(v,v)>0$ holds for all lightlike $v\in TM$ (a lightlike vector is nonzero by convention). For the following definition, note that when $M$ is a manifold-with-boundary and $x\in\mfbd M$, then each vector in $T_xM$ is either tangential to $\mfbd M$ or \emph{inward-directed} or \emph{outward-directed} in a well-defined sense.

\begin{definition}
Let $S_0,S_1$ be $(n-1)$-dimensional closed manifolds. A \emph{weak Lorentz cobordism between $S_0$ and $S_1$} is a compact $n$-dimensional Lorentzian man\-ifold-with-boundary $(M,g)$ whose boundary is the disjoint union $S_0\sqcup S_1$, such that $M$ admits a $g$-timelike vector field which is inward-directed on $S_0$ and outward-directed on $S_1$. A \emph{Lorentz cobordism between $S_0$ and $S_1$} is a weak Lorentz cobordism $(M,g)$ between $S_0$ and $S_1$ such that $\mfbd M$ is $g$-spacelike. $S_0$ is [\emph{weakly}] \emph{Lorentz cobordant} to $S_1$ iff there exists a [weak] Lorentz cobordism between $S_0$ and $S_1$.
\end{definition}

[Weak] Lorentz cobordance is an equivalence relation. Two manifolds are Lorentz cobordant if and only if they are weakly Lorentz cobordant. But when we require the cobordism to satisfy in addition the strict lightlike convergence condition, we obtain two extremely different cobordance relations.

\smallskip
A theorem of F.\ Tipler \cite{Tipler1977} which we review in Section \ref{tiplersection} implies that whenever two manifolds $S_0,S_1$ are Lorentz cobordant via a cobordism that satisfies the strict lightlike convergence condition, then $S_0,S_1$ are diffeomorphic. The situation is completely different for weak Lorentz cobordance:

\begin{theorem} \label{maincob}
Let $n\geq4$, let $S_0,S_1$ be closed $(n-1)$-dimensional manifolds, let $(M,g)$ be a weak Lorentz cobordism between $S_0$ and $S_1$, let $\Lambda\in\R$. If $n=4$, assume that $M$ is orientable and has no closed connected component. Then there exists a weak Lorentz cobordism $(M,g')$ between $S_0$ and $S_1$ such that
\begin{itemize}
\item
every $g$-causal vector in $TM$ is $g'$-timelike;
\item
$(M,g')$ satisfies the strict lightlike convergence condition and the dominant energy condition with respect to $\Lambda$;
\item
$M$ does not admit any codimension-one foliation none of whose tangent vectors is $g'$-timelike; in particular, $(M,g')$ does not admit any space foliation.
\end{itemize}
\end{theorem}

Again there is a weaker version for $3$-manifolds: Theorem \ref{maincob3} below.

\smallskip
Theorem \ref{maincob} implies in particular that for all orientable closed $3$-manifolds $S_0,S_1$, there exists an orientable weak Lorentz cobordism from $S_0$ to $S_1$ which satisfies the strict lightlike convergence condition. The contrast to Tipler's theorem is evident.

\subsection*{Acknowledgements}

The present article arose from a question that I had been asked by Christian Bär. I thank Felix Finster for several helpful comments, and I am grateful to Kai Zehmisch for a remark which eventually led me to reference \cite{Varela1976}.

\section{Preliminaries: Integrable and nonintegrable distributions} \label{backgroundI}

\begin{generalconventions}
$\N=\set{0,1,2,\dots}$. A \emph{manifold} does not have a boundary; a \emph{manifold-with-boundary} might have an empty boundary. All manifolds, bundles, sections in bundles etc.\ are assumed to be smooth, except when stated otherwise. All vector spaces/bundles are over the field $\R$. When $V,W$ are two such spaces/bundles, then $\Lin(V,W)$ denotes the vector space/bundle of linear maps $V\to W$. We use the terms \emph{semi-Riemannian/Lorentz[ian] metric, index, timelike, lightlike} as in \cite{ONeill}. We say that a vector is \emph{spacelike} iff it is spacelike in the sense of \cite{ONeill} and nonzero. A vector is \emph{causal} iff it is timelike or lightlike. A \emph{pseudo-Riemannian metric} on an $n$-manifold is a semi-Riemannian metric of index $\in\set{1,\dots,n-1}$.
\end{generalconventions}

\begin{definition}
Let $M$ be an $n$-manifold, let $p\in\set{0,\dots,n}$. We denote the fib{\re} over $x\in M$ of a $p$-plane distribution $H$ on $M$ by $H_x$ or $H(x)$. A \emph{line distribution} is a $1$-plane distribution. A $p$-plane distribution $H$ on $M$ is \emph{spacelike} [resp.\ \emph{timelike}] with respect to a semi-Riemannian metric $g$ on $M$ iff every nonzero vector in $H$ is spacelike [resp.\ timelike]. When $g$ has index $q$, then a \emph{space} [resp.\ \emph{time}] \emph{distribution} on $(M,g)$ is a spacelike [timelike] $(n-q)$-plane [resp.\ $q$-plane] distribution on $M$. A foliation on $(M,g)$ is called a \emph{spacelike} [resp.\ \emph{space}, \emph{timelike}, \emph{time}] \emph{foliation} iff its tangent distribution is a spacelike [space, timelike, time] distribution. (The fib{\re} over $x\in M$ of the \emph{tangent distribution} to a given foliation is the tangent space in $x$ to the leaf through $x$.) A distribution is \emph{integrable} iff it is the tangent distribution of a foliation. Two distributions $V,H$ on $M$ are \emph{complementary} iff the tangent bundle $TM$ is the internal direct sum of $V$ and $H$. The \emph{orthogonal distribution} of a $p$-plane distribution $H$ on $M$ with respect to a semi-Riemannian metric $g$ on $M$ is the $(n-p)$-plane distribution $\bot_gH$ on $M$ whose fib{\re} over $x$ is $\set{v\in T_xM \suchthat \forall w\in H_x\colon g(v,w)=0}$. When the orthogonal distribution $V$ of $H$ is complementary to $V$, then we call $V$ also the \emph{orthogonal complement} of $H$. (This happens for instance when $H$ is spacelike or timelike. When $H$ is a space [resp.\ time] distribution, then its orthogonal complement is a time [space] distribution.)
\end{definition}

\begin{definition}[twistedness]
Let $H$ be a $p$-plane distribution on an $n$-manifold $M$. The \emph{twistedness\footnote{As far as I can tell, the name \emph{twistedness} is due to W.\ Thurston: \cite{Thurston}, p.~176.} of $H$} is a section $\Twist_H$ in the vector bundle $\Lambda^2(H^\ast)\otimes(TM/H)$ (i.e., it is a $TM/H$-valued $2$-form on $M$), defined as follows: Let $\pi\colon TM\to TM/H$ denote the obvious projection. For all $x\in M$ and $v_0,v_1\in H_x$, we define
\[
\Twist_H(v,w) = \pi([\hat{v}_0,\hat{v}_1])(x) \;\;;
\]
here $\hat{v}_i$ is any section in $H$ with $\hat{v}_i(x)=v_i$, and $[.,.]$ denotes the Lie bracket of vector fields on $M$, so $\pi([\hat{v}_0,\hat{v}_1])$ is a section in $TM/H$. (Note that $\Twist_H(v,w)$ is well-defined, i.e.\ independent of the choice of $\hat{v}_0,\hat{v}_1$, because the map $(\hat{v}_0,\hat{v}_1)\mapsto \pi([\hat{v}_0,\hat{v}_1])$ is $C^\infty(M,\R)$-bilinear.) In a context where a complementary distribution $V$ of $H$ is specified, we will usually identify $TM/H$ with $V$, and thereby $\Twist_H$ with a section in $\Lambda^2(H^\ast)\otimes V$. When a semi-Riemannian metric is specified which makes $H$ spacelike or timelike, then we identify $TM/H$ with $\bot_gH$.

\smallskip
We call $H$ \emph{twisted at $x\in M$} iff the section $\Twist_H$ does not vanish at $x$, i.e., iff there exist $v,w\in H_x$ with $\Twist_H(v,w)\neq0$. We call $H$ \emph{twisted} iff it is twisted at every point of $M$.
\end{definition}

By the Frobenius theorem (\cite{Conlon}, Theorem 4.5.5), a distribution is integrable if and only if its twistedness vanishes identically. Thus $\Twist_H$ measures how far $H$ is from being integrable. Every line distribution $H$ is integrable because the vector bundle $\Lambda^2(H^\ast)$ has rank $0$.

\begin{distributionspace} \label{C0topdef1}
When $E\to M$ is a fib{\re} bundle, let $C^\infty(M\ot E)$ denote the set of smooth sections in $E$. This set can be equipped with several interesting topologies, but the most important one in the present article is the \emph{compact-open $C^0$-topology} (synonymously: \emph{compact-open topology} or \emph{topology of locally uniform convergence}), which we will simply call the \emph{$C^0$-topology} from now on. (The set $C^\infty(M\ot E)$ is a subset of the set $C^0(M,E)$ of all continuous maps [not necessarily sections] from $M$ to $E$. The $C^0$-topology on $C^\infty(M\ot E)$ is just the subspace topology induced by the usual compact-open topology on $C^0(M,E)$.)

\smallskip
When $M$ is an $n$-manifold and $p\in\set{0,\dots,n}$, we can take $E\to M$ to be the \emph{Grassmann bundle} $\Gr_p(TM)\to M$, whose fib{\re} over $x$ is the Grassmann manifold $\Gr_p(T_xM)$, i.e.\ the ($p(n-p)$-dimensional) space of all $p$-dimensional sub vector spaces of $T_xM$. A $p$-plane distribution on $M$ is just a section in $\Gr_p(TM)\to M$. We equip the set $\Distr_p(M)\define C^\infty(M\ot\Gr_p(TM))$ of $p$-plane distributions on $M$ with the $C^0$-topology.
\end{distributionspace}

\begin{notation} \label{C0notation}
Let $V,W$ be complementary distributions on a manifold $M$.

\smallskip
For every distribution $Z$ on $M$ which is complementary to $V$, let $\lambda[Z]\colon W\to V$ denote the vector bundle morphism given as follows: For every $x\in M$, the restriction $\lambda[Z]_x$ of $\lambda[Z]$ to the fib{\re} $W_x$ is the unique linear map $W_x\to V_x$ whose graph $\subseteq W_x\times V_x = T_xM$ is $Z_x$. I.e., for all $w\in W_x$, the vector $\lambda[Z](w)$ is the unique $v\in V_x$ with $w+v\in Z_x$.
\end{notation}

\begin{distributionspace}[again] \label{C0topdef2}
For the proofs in the sections \ref{C0closed} and \ref{main0proof}, we need a more explicit description of the $C^0$-topology on $\Distr_p(M)$. Since this topology is metri{\z}able, it suffices to say when precisely a sequence $(H_k)_{k\in\N}$ in $\Distr_p(M)$ converges in $C^0$ to a distribution $H\in\Distr_p(M)$. (Even if the topology were not metri{\z}able, this information is all we need below.)

\smallskip
Let $V,W$ be complementary distributions on $M$ such that $V$ is complementary to $H$. We use the notation $\lambda[.]$ from \ref{C0notation} with respect to these data.

\smallskip
The sequence $(H_k)_{k\in\N}$ converges to $H$ with respect to the $C^0$-topology if and only if the following conditions hold for every compact subset $K$ of $M$:
\begin{itemize}
\item
there exists a number $k_K\in\N$ such that for all $k\geq k_K$, the distribution $H_k$ is complementary to $V$ on the set $K$;
\item
$\displaystyle{\lim_{k_K\leq k\to\infty}}\norm{\lambda[H_k]-\lambda[H]}_{C^0(K,\Lin(W,V))} = 0$.
\end{itemize}
(In order to define $\norm{.}_{C^0(K,\Lin(W,V))}$, we should choose a Riemannian metric $h$ on $M$. This induces fib{\re}wise norms on the vector bundles $W,V$ and thus a fib{\re}wise operator norm on $\Lin(W,V)$. But since $K$ is compact, all $h$ yield equivalent operator norms and thus the same convergence criterion.)

\smallskip
It is not hard to see that this convergence criterion does not depend on the choice of $V,W$. We leave it to the reader to check carefully that the $C^0$-topology defined in \ref{C0topdef1} is really the same as the one described here (by spelling out how the topology of the Grassmann manifold is defined and how the fib{\re} bundle structure of the Grassmann bundle is induced by the vector bundle structure of $TM$). Our results below do not depend on this fact because they employ only the definition given here in \ref{C0topdef2}.
\end{distributionspace}

We will repeatedly use the following basic fact (cf.\ e.g.\ \cite{Baum}, Satz 0.48):

\begin{theorem}
Every semi-Riemannian manifold $(M,g)$ admits a time distribution $V$ and thus also a space distribution (e.g.\ the orthogonal complement of $V$).
\end{theorem}

More generally, every time distribution on a closed subset $A$ of $M$ can be extended to a time distribution on $M$, but we do not need that in the present article.

\begin{conncomp} \label{connectedcomponents}
Let $M$ be an $n$-manifold. Considering semi-Riemannian metrics on $M$ as sections in the vector bundle of symmetric bilinear forms on $TM$, we can equip the set $\Metr_q(M)$ of all index-$q$ metrics on $M$ with the $C^0$-topology. The resulting topological space is locally path-connected. The set $\pi_0(\Metr_q(M))$ of its connected components is in canonical bijective correspondence to the set of connected components of $\Distr_q(M)$: For each connected component $\mathscr{C}$ of $\Metr_q(M)$, we choose a metric $g$ in $\mathscr{C}$ and a $g$-time distribution $V\in\Distr_q(M)$, and we assign to $\mathscr{C}$ the connected component of $V$ in $\Distr_q(M)$. This yields a well-defined bijection. Its inverse is obtained by choosing to each connected component $\mathscr{C}'$ of $\Distr_q(M)$ a distribution $V$ in $\Distr_q(M)$ and an index-$q$ metric $g$ which makes $V$ timelike, and assigning to $\mathscr{C}'$ the connected component of $g$ in $\Metr_q(M)$. Similarly, we obtain a canonical bijection between $\pi_0(\Metr_q(M))$ and $\pi_0(\Distr_{n-q}(M))$ by replacing timelike distributions with spacelike distributions in the description above. (Details can be found in \cite{Nardmann2004}, Appendix D.)

\smallskip
In particular, in the situation of Theorems \ref{main0}, \ref{maindom}, \ref{maincob}, the metrics $g$ and $g'$ lie in the same connected component of $\Metr_q$ resp.\ $\Metr_1(M)$, because there exists a distribution of rank $q$ resp.\ $1$ which they both make timelike.
\end{conncomp}

\section{Preliminaries: Energy conditions} \label{backgroundE}

\begin{definition}
Let $(M,g)$ be a semi-Riemannian manifold, let $\Lambda\in\R$. The \emph{energy-momentum tensor of $(M,g)$ with respect to (the cosmological constant) $\Lambda$} is the symmetric bilinear form field $T\define \Ric_g -\frac{1}{2}\scal_gg +\Lambda g$.

\smallskip
\emph{$(M,g)$ satisfies the weak energy condition with respect to $\Lambda$} iff $T(v,v)\geq0$ holds for every $g$-timelike vector $v\in TM$. \emph{$(M,g)$ satisfies the semi-dominant energy condition with respect to $\Lambda$} iff for every $g$-timelike vector $v\in TM$, the vector $-\sharp(T(v,.))$ is not spacelike; here $\sharp\colon T^\ast M\to TM$ denotes the isomorphism induced by $g$. \emph{$(M,g)$ satisfies the dominant energy condition with respect to $\Lambda$} iff it satisfies the weak energy condition and the semi-dominant energy condition with respect to $\Lambda$.

\smallskip
(If you are a physicist accustomed to a certain unit system, you might prefer the energy-momentum tensor to be defined by $cT = \Ric_g -\frac{1}{2}\scal_gg +\Lambda g$, where $c>0$ is a constant depending on the unit system; e.g.\ $c=8\pi$. Such a constant does not change the energy conditions and is therefore irrelevant here. Note that the term \emph{semi-dominant energy condition} is my invention; the concept does not seem to have a standard name.)
\end{definition}

\begin{remark}
When $g$ is Lorentzian, the dominant energy condition can be stated in an obviously equivalent way: For every $x\in M$, the set of timelike vectors in $T_xM$ has two connected components. The dominant energy condition holds iff for every $g$-timelike vector $v\in TM$, the vector $-\sharp(T(v,.))$ lies in the closure of the connected component which contains $v$.
\end{remark}

In General Relativity, where $g$ is Lorentzian, these conditions have a clear physical motivation: Consider an observer who moves through the point $x\in M$ with $4$-velocity $v\in T_xM$ (which by definition satisfies $g(v,v)=-1$). Then the number $T(v,v)$ is the mass (i.e.\ energy) density of matter in the point $x$. The weak energy condition says roughly that mass should be nonnegative. When $z\in T_pM$ is orthogonal to $v$ with $g(z,z)=1$, then $-T(v,z)$ is the $z$-component of the momentum density of matter in the point $x$ as seen by the observer. Thus $-\sharp(T(v,.))$ is the $4$-momentum (i.e.\ energy-momentum) density of matter in the point $x$ as seen by the observer. The semi-dominant energy condition says roughly that matter does not move faster than light.

\begin{definition} \label{lccdef}
A semi-Riemannian manifold $(M,g)$ satisfies the [\emph{strict}] \emph{lightlike convergence condition} iff $\Ric_g(v,v)\geq0$ [resp.\ $\Ric_g(v,v)>0$] holds for every lightlike vector $v\in TM$. It satisfies the [\emph{strict}] \emph{timelike convergence condition} iff $\Ric_g(v,v)\geq0$ [resp.\ $\Ric_g(v,v)>0$] holds for every timelike vector $v\in TM$. It satisfies the \emph{strict causal convergence condition} iff it satisfies the strict lightlike convergence condition and the strict timelike convergence condition.
\end{definition}

When $(M,g)$ satisfies the weak energy condition with respect to some $\Lambda$, or when it satisfies the timelike convergence condition, then it satisfies the lightlike convergence condition.

\smallskip
Let us also introduce strict versions of the weak and semi-dominant energy conditions:

\begin{definition}
Let $(M,g)$ be a semi-Riemannian manifold, let $\Lambda\in\R$, let $T$ denote the energy-momentum tensor of $(M,g)$ with respect to $\Lambda$. \emph{$(M,g)$ satisfies the strict weak energy condition with respect to $\Lambda$} iff $T(v,v)>0$ holds for every $g$-causal vector $v\in TM$. \emph{$(M,g)$ satisfies the strict semi-dominant energy condition with respect to $\Lambda$} iff for every $g$-causal vector $v\in TM$, the vector $-\sharp(T(v,.))$ is timelike. \emph{$(M,g)$ satisfies the strict dominant energy condition with respect to $\Lambda$} iff it satisfies the strict weak energy condition and the strict semi-dominant energy condition with respect to $\Lambda$.
\end{definition}

When $(M,g)$ satisfies the strict weak energy condition with respect to some $\Lambda$, then it satisfies the strict lightlike convergence condition.

\section{The integrability property is $C^0$-closed} \label{C0closed}

The main tool in the proof of Theorem \ref{main0} is the following proposition. It would remain true if we considered $C^1$ distributions instead of smooth distributions (actually, even less regularity would suffice), but that is not important for our purposes.

\begin{proposition} \label{propclosed}
Let $p\in\N$. Whenever a sequence of (smooth) integrable $p$-plane distributions on a given manifold converges to a (smooth) $p$-plane distribution $H$ with respect to the $C^0$-topology, then $H$ is integrable, too.
\end{proposition}

Note that the proposition becomes trivial when we replace the $C^0$-topology by the $C^1$-topology, because the function which maps each distribution on a compact manifold to the $C^0$-norm of its twistedness is continuous with respect to the $C^1$-topology (thus vanishing twistedness of each element in the sequence implies vanishing twistedness of the limit). It is not continuous with respect to the $C^0$-topology, so the proposition might be surprising at first sight.

\smallskip
As far as I know, this rather elementary result is not mentioned explicitly in the literature. The case of codimension-one distributions appears in an article of F.\ Varela (\cite{Varela1976}, p.~255; note that a nonintegrable codimension-one distribution is just the kernel of a $1$-form $\omega$ with $\omega\wedge d\omega\not\equiv0$), but with only a rough sketch of proof. Varela cites (5.2 on p.~242) an old paper of G.\ Reeb as a reference, but that article seems to contain only a vague remark related to the issue. If I interpret Varela's sketch correctly, his argument employs Darboux' theorem on contact forms in dimension $3$. Since Darboux' theorem does not generali{\z}e from contact structures to nonintegrable distributions of higher codimension, neither does Varela's argument.

\smallskip
It seems therefore appropriate to give a detailed proof of Proposition \ref{propclosed}. That is basically an exercise in the theory of ordinary differential equations. At the end of the section, we review briefly Varela's argument and compare it to the codimension-one special case of the proof given here.

\begin{proof}[Proof of Proposition \ref{propclosed}]
Let $(H_k)_{k\in\N}$ be a sequence of integrable $p$-plane distributions on an $n$-manifold $M$ which converges locally uniformly to a distribution $H$, and let $x\in M$. We have to show that $H$ is integrable on some open neighb{\ou}rhood $U$ of $x$. We can assume $M=\R^n$ and $x=0$ and $H(x)=\R^p\times\set{0}\subseteq \R^p\times\R^q = T_x\R^n$ without loss of generality. Let $q\define n-p$, and let $V$ denote the $q$-plane distribution on $M$ given by $V(x)=\set{0}\times\R^q\subseteq T_x\R^n$ for all $x$.

\medskip
For all but finitely many $k\in\N$, the subspace $H_k(x)$ is transverse to $V(x)$ (because the $H_k(x)$ converge to $H(x)$). By deleting the finitely many other $k$ from the sequence, we arrange that $H_k(x)$ is transverse to $V(x)$ for all $k$. For each $k$, we consider the set $U_k$ of all $y\in M$ such that $H_k(y)$ is complementary to $V(y)$; clearly $U_k$ is a neighb{\ou}rhood of $x$. The intersection of all $U_k$ contains a compact neighb{\ou}rhood $K$ of $x$. (Otherwise there would exist a sequence $(x_j)_{j\in\N}$ in $M$ converging to $x$ and a sequence $(k_j)_{j\in\N}$ in $\N$ such that $H_{k_j}(x_j)\cap V(x_j)$ contains a vector $v_j\in T_{x_j}\R^n=\R^n$ of euclidean norm $1$. The sequence $(k_j)_{j\in\N}$ contains a subsequence $(l_j)_{j\in\N}$ converging to $\infty$, for otherwise a number $k$ would occur infinitely often in $(k_j)_{j\in\N}$, in contradiction to $U_k$ being a neighb{\ou}rhood of $x$. By deleting the other elements from the sequence $(k_j)_{j\in\N}$, we may assume $k_j=l_j$ for all $j$. Since the distributions $H_{k_j}$ converge locally uniformly to $H$, the $p$-planes $H_{k_j}(x_j)$ converge to $H(x)$. The sequence $(v_j)_{j\in\N}$ in the compact set $S^{n-1}\cap(\set{0}\times\R^q)\subseteq\R^n$ has a subsequence which converges to some nonzero vector $v\in V(x)$. Since $v_j\in H_{k_j}(x_j)$, we obtain $v\in H(x)\cap V(x) = \set{0}$, a contradiction.)

\smallskip
In the following, the notations $\lambda[.]$ (cf.\ \ref{C0notation}) and $\norm{.}_{C^0(K,\Lin(W,V))}$ (which we abbreviate as $\norm{.}_{C^0(K)}$ or $\norm{.}_{C^0}$) refer to the decomposition $W\oplus V$ of $TM\restrict K$ which is fib{\re}wise the decomposition $T_x\R^n=\R^p\oplus\R^q$, and to the euclidean metric $h$ on $M=\R^n$ (cf.\ the remark in \ref{C0topdef2}).

\smallskip
The uniform convergence $H_k\to H$ on $K$ implies that $a\define\sup\set{\norm{\lambda[H_k]}_{C^0(K)} \suchthat k\in\N}$ $\in[0,\infty]$ is finite, and that $\norm{\lambda[H]}_{C^0(K)}\leq a$.

\smallskip
We choose a neighb{\ou}rhood $B\times J\subseteq K$ of $x=(0,0)\in\R^p\times\R^q$, where $B\subseteq\R^p$ is a closed ball of radius $r_B$ centered at $0\in\R^p$, and $J\subseteq\R^q$ is a closed ball of radius $r_J$ centered at $0\in\R^q$. Since $H$ is Lipschitz continuous on $B\times J$, the $\Lin(\R^p,\R^q)$-valued function $\lambda[H]$ is Lipschitz continuous on $B\times J$ as well, i.e., there exists a number $C>0$ such that $\abs{\lambda[H(y_0)]-\lambda[H(y_1)]} \leq C\abs{y_0-y_1}$ for all $y_0,y_1\in B\times J$. By shrinking $B$ if necessary, we arrange that $r_BC < 1$ and $ar_B < r_J$. Let $r_I$ be the number $r_J-ar_B>0$, and let $I\subseteq J$ be the closed $0$-centered ball of radius $r_I$ in $\R^q$.

\smallskip
For each $k\in\N$ and $t\in I$, the intersection of $B\times\R^q$ and the leaf through the point $(0,t)\in B\times I$ of the integrable distribution $H_k$ is the graph of a unique function $f^k_t\in C^\infty(B,J)$. (Since $H_k$ is transverse to $V$ on $B\times J$, there are only two alternatives: either the intersection of $B\times\R^q$ and the leaf through $(0,t)$ is the graph of a function $\in C^\infty(B,J)$ for each $t\in I$; or there is a $t\in I$ such that the intersection of $\mathring{B}\times\R^q$ and the leaf through $(0,t)$ meets $B\times\mfbd J$. The second alternative is ruled out by $r_J-\abs{t} \geq r_J-r_I = ar_B$: Let $\gamma=(\gamma_B,\gamma_J)\colon[0,1]\to B\times J$ be a smooth path from $(0,t)$ to a point $(z,\tau)\in \mathring{B}\times\mfbd J$, i.e.\ a point $(z,\tau)$ with $\abs{z}< r_B$ and $\abs{\tau}=r_J$, such that $\im(\gamma)$ is contained in the leaf through $(0,t)$. Since $\gamma'(s)\in H_k(\gamma(s))$ and thus $\abs{\gamma_J'(s)}\leq \norm{\lambda[H_k]}_{C^0(K)}\,\abs{\gamma_B'(s)} \leq a\abs{\gamma_B'(s)}$ holds for all $s\in[0,1]$, we obtain
\[
r_J-\abs{t} \leq \abs{\gamma_J(1)-t} \leq \int_0^1\abs{\gamma_J'(s)}ds \leq a\int_0^1\abs{\gamma_B'(s)}ds \leq a(\gamma_B(1)-\gamma_B(0)) < ar_B \;,
\]
a contradiction.)

\smallskip
By definition, the map $f^k_t$ satisfies $f^k_t(0)=t$, and the subspace $H_k(z,f^k_t(z))\subseteq\R^n$ is the graph of the derivative $D_zf^k_t\in \Lin(\R^p,\R^q)$; i.e., $D_zf^k_t = \lambda[H_k(z,f^k_t(z))]$.

\smallskip
We claim that the $f^k_t$ converge uniformly to a function $f_t\in C^0(B,J)$ as $k$ tends to infinity. In order to prove this, we define, for all $k\in\N$ and $(z,t)\in B\times I$, the function $\alpha^k_{z,t}\colon [0,1]\to J$ by $\alpha^k_{z,t}(s) = f^k_t(sz)$. Then $f^k_t(z)=\alpha^k_{z,t}(1)$, and $\alpha^k_{z,t}$ solves the ordinary differential equation
\[ \begin{split}
(\alpha^k_{z,t})'(s) &= (D_{sz}f^k_t)(z) = \lambda[H_k(sz,f^k_t(sz))](z) = \lambda[H_k(sz,\alpha^k_{z,t}(s))](z)
\end{split} \]
with initial value $\alpha^k_{z,t}(0) = f^k_t(0) = t$.

\smallskip
Let $k,l\in\N$ and $(z,t)\in B\times I$. The map $[0,1]\to\R$ given by $s\mapsto \abs{\alpha^k_{sz,t}(1)-\alpha^l_{sz,t}(1)}$ is continuous and thus achieves its maximum in some $\bar{s}\in[0,1]$. For $\bar{z}\define\bar{s}z$, we compute
\[ \begin{split}
&\abs{\alpha^k_{\bar{z},t}(1)-\alpha^l_{\bar{z},t}(1)}
\leq \int_0^1\Big|{\lambda[H_k(s\bar{z},\alpha^k_{\bar{z},t}(s))](\bar{z}) -\lambda[H_l(s\bar{z},\alpha^l_{\bar{z},t}(s))](\bar{z})}\Big|ds\\
&\leq \int_0^1\Big|{\lambda[H_k(s\bar{z},\alpha^k_{\bar{z},t}(s))](\bar{z}) -\lambda[H(s\bar{z},\alpha^k_{\bar{z},t}(s))](\bar{z})}\Big|ds\\
&\mspace{20mu}+\int_0^1\Big|{\lambda[H_l(s\bar{z},\alpha^l_{\bar{z},t}(s))](\bar{z}) -\lambda[H(s\bar{z},\alpha^l_{\bar{z},t}(s))](\bar{z})}\Big|ds\\
&\mspace{20mu}+\int_0^1\Big|{\lambda[H(s\bar{z},\alpha^k_{\bar{z},t}(s))](\bar{z}) -\lambda[H(s\bar{z},\alpha^l_{\bar{z},t}(s))](\bar{z})}\Big|ds\\
&\leq \abs{\bar{z}}\bigg(\norm{\lambda[H_k]-\lambda[H]}_{C^0} +\norm{\lambda[H_l]-\lambda[H]}_{C^0} +C\int_0^1\Big|{\alpha^k_{\bar{z},t}(s)-\alpha^l_{\bar{z},t}(s)}\Big|ds\bigg)\\
&= \abs{\bar{z}}\bigg(\norm{\lambda[H_k]-\lambda[H]}_{C^0} +\norm{\lambda[H_l]-\lambda[H]}_{C^0}\!\bigg) +\abs{\bar{z}}C\!\!\int_0^1\Big|{\alpha^k_{s\bar{z},t}(1)-\alpha^l_{s\bar{z},t}(1)}\Big|ds\\
&\leq r_B\bigg(\norm{\lambda[H_k]-\lambda[H]}_{C^0} +\norm{\lambda[H_l]-\lambda[H]}_{C^0}\bigg) +r_BC\Big|{\alpha^k_{\bar{z},t}(1)-\alpha^l_{\bar{z},t}(1)}\Big| \;\;;
\end{split} \]
here $\norm{.}_{C^0}$ is the $C^0$-norm for $\Lin(\R^p,\R^q)$-valued maps on the compact set $B\times I$. In the last step, we used $\abs{\alpha^k_{s\bar{z},t}(1)-\alpha^l_{s\bar{z},t}(1)}
= \abs{\alpha^k_{s\bar{s}z,t}(1)-\alpha^l_{s\bar{s}z,t}(1)}
\leq\abs{\alpha^k_{\bar{z},t}(1)-\alpha^l_{\bar{z},t}(1)}$.

\smallskip
We obtain
\[ \begin{split}
\abs{f^k_t(z)-f^l_t(z)} &= \abs{\alpha^k_{z,t}(1)-\alpha^l_{z,t}(1)}
\leq \abs{\alpha^k_{\bar{z},t}(1)-\alpha^l_{\bar{z},t}(1)}\\
&\leq \frac{r_B}{1-r_BC}\bigg(\norm{\lambda[H_k]-\lambda[H]}_{C^0} +\norm{\lambda[H_l]-\lambda[H]}_{C^0}\bigg) \;\;.
\end{split} \]

\smallskip
Since $(\lambda[H_k])_{k\in\N}$ converges uniformly on $B\times I$ to $\lambda[H]$, we see that $(f^k_t)_{k\in\N}$ is a Cauchy sequence in the Banach space $C^0(B,\R)$ and thus has a limit $f_t\in C^0(B,\R)$. It lies in $C^0(B,J)$ because $J$ is closed. This proves our claim from above.

\smallskip
Now we claim that the sequence $(Df^k_t)_{k\in\N}$ in $C^\infty(B,\Lin(\R^p,\R^q))$ converges uniformly to $\eta_t\define \big(z\mapsto\lambda[H(z,f_t(z))]\big)$. In fact, this follows immediately from the following inequality which holds for all $z\in B$:
\[ \begin{split}
&\abs{D_zf^k_t -\lambda[H(z,f_t(z))]}\\
&= \abs{\lambda[H_k(z,f^k_t(z))] -\lambda[H(z,f_t(z))]}\\
&\leq \abs{\lambda[H_k(z,f^k_t(z))] -\lambda[H(z,f^k_t(z))]} +\abs{\lambda[H(z,f^k_t(z))] -\lambda[H(z,f_t(z))]}\\
&\leq \norm{\lambda[H_k] -\lambda[H]}_{C^0} +C\abs{f^k_t(z)-f_t(z)} \;\;.
\end{split} \]

\smallskip
So for each $t\in I$, the derivatives of the $f^k_t$ converge uniformly to $\eta_t$, and the $f^k_t$ converge uniformly to $f_t$. Hence $f_t$ is $C^1$, the sequence $(f^k_t)_{k\in\N}$ converges in $C^1$ to $f_t$, and $Df_t = \eta_t = \lambda[H(.,f_t(.))]$. Thus the graph of $f_t$ is an integral manifold for the smooth distribution $H$. In particular, the Lie bracket of two local sections in $H$ is always a local section in $H$. By the Frobenius theorem, $H$ is thus integrable on $B\times I$, i.e.\ on a neighb{\ou}rhood of $x$.
\end{proof}

As announced above, we will now briefly explain Varela's alternative proof in the codimen\-sion-one case. It works by contradiction: We assume that there exists a sequence of integrable $(n-1)$-plane distributions on $\R^n$ which converges in $C^0$ to a nonintegrable distribution.

\smallskip
The first step is a reduction to the case $n=3$: If such a sequence exists for some $n$, then also for $n=3$. The argument is that each codimension-one distribution $H$ is the kernel of a nowhere vanishing $1$-form $\omega$, and $H$ is integrable if and only if $\omega\wedge d\omega=0$. Our limit distribution on $\R^n$ is nonintegrable, so there exist $x\in\R^n$ and $u,v,w\in T_x\R^n$ with $(\omega\wedge d\omega)(u,v,w)\neq0$. We restrict all $1$-forms to the $3$-plane spanned by $u,v,w$ and thereby get the desired sequence for $\R^3$.

\smallskip
Now Darboux' theorem on contact structures (cf.\ e.g.\ \cite{Geigessurvey}, Theorem 2.24) says that there exist local coordinates $(x_0,x_1,x_2)$ on some open set in $\R^3$ (without loss of generality on $\R^3$ itself) such that our limit distribution is $H=\ker(dx_0+x_2dx_1)$ there.

\smallskip
We define $Y\colon\R^3\to\R^2$ by $(z_0,z_1,z_2)\mapsto (-z_2,z_1)$. Let $\pr_{12}\colon\R^3\to\R^2$ denote the projection to the last two components. There is a unique vector field $X$ on $\R^3$ which is a section in $H$ and satisfies $\pr_{12}\compose X=Y$. In the same way, we define for each distribution $H_k$ in our $C^0$-converging sequence a vector field $X_k$. (For transversality reasons, the latter definition works only on some neighb{\ou}rhood of $0\in\R^3$. That the neighb{\ou}rhood can be chosen independent of $k$ is shown in the same way as in our proof above.)

\smallskip
The $\pr_{12}$-image of each integral curve of $X$ or $X_k$ is a (possibly degenerated) circle in $\R^2$ with cent{\re} $0$. The integral curves of $X$ can be computed explicitly (that's why we reduced to dimension $3$ and applied Darboux' theorem in the first place); except in degenerate cases, they are not closed but spiral-shaped.

\smallskip
On the other hand, each integral curve of each $X_k$ is closed because it stays within one leaf of the $2$-dimensional foliation defined by $H_k$.

\smallskip
Whenever a sequence $(X_k)_{k\in\N}$ of vector fields on $\R^n$ converges in $C^0$ to a (locally Lip\-schitz) vector field $X$ (and that's what's happening in our case), then the $X_k$-integral curves with initial value $x$ converge in $C^0$ to the $X$-integral curve with initial value $x$. But a sequence of closed curves cannot $C^0$-converge to a nonclosed curve, so we get a contradiction.\qed

\medskip
When one carries out all details of this proof, it is not shorter than the one we gave above (although Darboux' theorem is just cited here). In fact, several technical points arise in both proofs similarly.

\section{Proof of Theorem \ref{main0}} \label{main0proof}

\begin{lemma} \label{lemma0}
Let $(M,g)$ be an $n$-dimensional pseudo-Riemannian manifold of index $q\in\set{1,\dots,n-2}$, let $x\in M$. Then there exists a $g$-space distribution $H$ on $M$ which is twisted at $x$.
\end{lemma}
\begin{proof}
We choose a $g$-space distribution $H'$ on $M$. There exists an open neighb{\ou}rhood $U$ of $x$ on which $TM$ admits a frame $(e_1,\dots,e_n)$ such that $(e_1,\dots,e_{n-q})$ is a local frame of $H'$. If the value of the Lie bracket $[e_1,e_2]$ in $x$ is not contained in $H'_x$, then $H'$ is twisted at $x$, and we can take $H=H'$. Otherwise we choose a compact neighb{\ou}rhood $U'\subset U$ of $x$ and a function $\beta\in C^\infty(M,\R)$ with $\beta\restrict(M\without U')=0$ and $\beta(x)=0$ and $d_x\beta(e_2)=1$. For $c\in\R_{>0}$, we define a frame $(e^c_1,\dots,e^c_n)$ of $TU$ by $e^c_1\define e_1-c\beta e_n$ and $e^c_i\define e_i$ for $i>0$. We define the $(n-q)$-plane distribution $H^c$ on $M$ to be $H'$ on $M\without U'$, and to be $\spann(e^c_1,\dots,e^c_{n-q})$ on $U$; note that both definitions agree on the overlap. For every sufficiently small $c>0$, this $H^c$ is $g$-spacelike. It is twisted at $x$ because $[e^c_1,e^c_2](x) = [e_1,e_2](x) +c\,e_n(x)$ is not contained in $H^c_x = H'_x$.
\end{proof}

As we will see later, the preceding lemma can be sharpened considerably via Thom's jet transversality theorem or Gromov's h-principles: one can find space distributions $H$ which are twisted on large subsets of $M$.

\begin{notation}
Let $V,H$ be complementary distributions on the manifold $M$, let $g_V,g_H$ be semi-Riemannian metrics on the vector bundles $V,H$, respectively. Then $g_V\oplus g_H$ denotes the semi-Riemannian metric on $M$ whose (pointwise) restriction to $V$ is $g_V$, whose restriction to $H$ is $g_H$, and which makes $V$ and $H$ orthogonal to each other.
\end{notation}

\begin{keylemma} \label{keylemma}
Let $M$ be a manifold whose tangent bundle has a splitting $TM=V\oplus H$ such that the $p$-plane distribution $H$ is twisted at $x\in M$. Let $g_V,g_H$ be Riemannian metrics on the vector bundles $V,H$, respectively. Let $U$ be an open neighb{\ou}rhood of $x$ in $M$. Let $(f_k)_{k\in\N}$ be a sequence of functions $f_k\in C^\infty(M,\R_{>0})$ such that $f_k\geq k$ on $U$. For $k\in\N$, let $g_k$ denote the pseudo-Riemannian metric $(-f_k^2\,g_V)\oplus g_H$ on $M$. Then there exists a $k_U\in\N$ such that for all $k\geq k_U$, the open set $U$ does not admit a $p$-dimensional foliation none of whose tangent vectors is $g_k$-timelike.
\end{keylemma}
\begin{proof}
Assume that no $k_U$ with the stated property exists. Then there is a sequence $(k(j))_{j\in\N}$ in $\N$ which converges to $\infty$, such that for all $j\in\N$, the set $U$ admits an integrable $p$-plane distribution $H_{k(j)}$ none of whose tangent vectors is $g_{k(j)}$-timelike.

\smallskip
We claim that the sequence $(H_{k(j)})_{j\in\N}$ converges to $H$ (actually, to the restriction of $H$ to $U$, but we suppress that in the notation) with respect to the $C^0$-topology on $\Distr_p(U)$. This follows from the convergence criterion given in \ref{C0topdef2}. Namely, $H_{k(j)}$ is complementary to $V$ for all $j$. We choose the Riemannian metric $h = g_H\oplus g_V$ in order to define the notations $\lambda[.]$ and $\norm{.}_K\define \norm{.}_{C^0(K,\Lin(W,V))}$.

\smallskip
By definition, we have $\abs{H_{k(j)}}(x) = \max\big\{ \abs{\lambda[H_{k(j)}](w)}_{g_V} \;\big|\; w\in H_x, \;\abs{w}_{g_H}\leq1 \big\}$, and $\lambda[H_{k(j)}](w)$ is the unique vector $v\in V_x$ with $v+w\in H_{k(j)}$. Since no vector in $H_{k(j)}$ is timelike with respect to $(-f_{k(j)}^2g_V)\oplus g_H$, we obtain $-f_{k(j)}^2g_V(v,v) +g_H(w,w)\geq0$ and thus $\abs{\lambda[H_{k(j)}](w)}_{g_V} \leq f_{k(j)}^{-1}\abs{w}_{g_H} \leq \frac{1}{k(j)}$. In particular, $\lim_{j\to\infty}\norm{H_{k(j)}}_K=0$ for every compact set $K$ in $U$. This proves our claim that $(H_{k(j)})_{j\in\N}$ converges in $C^0$ to $H$.

\smallskip
Now Proposition \ref{propclosed} implies that $H$ is integrable on $U$, in contradiction to the assumption that $H$ is twisted at $x$.
\end{proof}

The preceding lemma says, intuitively speaking, that space foliations cease to exist when one squeezes the spacelike region in the tangent bundle in such a way that it becomes concentrated around a nonintegrable distribution. This squeezing can also be roughly understood as increasing the speed of light.

\smallskip
The proof of the lemma shows why we had to prove the $C^0$-closedness of the integrability condition in Section \ref{C0closed}, i.e., why it would have been not enough to know the completely obvious $C^1$-closedness of this condition: The only information that we had about the distributions $H_k$ in the proof was their nontimelikeness, which yields only information about the $C^0$-topology.

\begin{proof}[\textit{Proof of Theorem \ref{main0}}]
By Lemma \ref{lemma0}, there exists a $g$-space distribution $H$ on $M$ which is twisted at some point $x\in M\without A$. Let $V$ denote the $g$-orthogonal distribution of $H$, let $g_V$ denote the Riemannian metric $-g\restrict V$ on the vector bundle $V$, and let $g_H$ denote the Riemannian metric $g\restrict H$ on $H$. We choose a function $\beta\in C^\infty(M,[0,1])$ with $\beta=0$ on $A$, and $\beta=1$ on a compact neighb{\ou}rhood $B\subset M\without A$ of $x$.

\smallskip
For every $k\in\N$, let $f_k$ be the function $\sqrt{1+k^2\beta^2}\in C^\infty(M,\R_{\geq1})$, and let $g_k$ denote the pseudo-Riemannian metric $(-f_k^2g_V)\oplus g_H$ of index $q$ on $M$. Clearly for all $k\in\N$, every $g$-timelike vector in $TM$ is $g_k$-timelike, and $g=g_k$ holds on $A$. By the Key lemma, there exists a $k\in\N$ such that the interior of $B$ does not admit a codimension-$q$ foliation none of whose tangent vectors is $g_k$-timelike.
\end{proof}

\section{The explicit example $g^c_n$ on $\R^n$} \label{specialcase}

We are now going to prove Remark \ref{mainremark} and Theorem \ref{mainRn}.

\begin{proof}[Proof of the first three statements in Remark \ref{mainremark}]
With respect to the standard coordinates $(x_0,x_1,x_2)$, the metric $g^c_3$ is given by the matrix-valued function
\[
(g_{ij})_{i,j\in\set{0,1,2}}\define \begin{pmatrix}
-c^2 &-c^2x_2 &0\\
-c^2x_2 &1-c^2x_2^2 &0\\
0 &0 &1
\end{pmatrix}
\]
because $\partial_0=ce_0$, $\partial_1=e_1+x_2ce_0$, $\partial_2=e_2$. Consider the diffeomorphism $\varphi\colon\R^n\to\R^n$ given by
\[
(x_0,x_1,x_2,\dots,x_{n-1})\mapsto (c^2x_0,cx_1,cx_2,\dots,cx_{n-1}) \;\;.
\]
The derivative $D_x\varphi\in\Lin(T_x\R^n,T_{\varphi(x)}\R^n)=\Lin(\R^n,\R^n)$ is given by the diagonal matrix $\diag(c^2,c,c,\dots,c)$. The metric $\tilde{g}\define\varphi^\ast(g^1_n)$ is determined by the values $\tilde{g}_x(\partial_i,\partial_j) = (g^1_n)_{\varphi(x)}((D_x\varphi)(\partial_i),(D_x\varphi)(\partial_j))$, i.e., its value in $x$ is given by the matrix
\begin{align*}
(\tilde{g}_{ij})_x &= \begin{pmatrix}c^2&&&\\ &c&&\\ &&c&\\ &&&\ddots \end{pmatrix}
\begin{pmatrix} -1 &-cx_2 &&\\ -cx_2 &1-c^2x_2^2 &&\\ &&1&\\ &&&\ddots
\end{pmatrix}
\begin{pmatrix}c^2&&&\\ &c&&\\ &&c&\\ &&&\ddots \end{pmatrix}\\
&= \begin{pmatrix} -c^4 &-c^4x_2 &&\\ -c^4x_2 &c^2(1-c^2x_2^2) &&\\ &&c^2&\\ &&&\ddots
\end{pmatrix}
= c^2((g^c_n)_{ij})_x \;\;.
\end{align*}
This proves $g^c_n = \frac{1}{c^2}\varphi^\ast(g^1_n)$, i.e.\ the first statement of Remark \ref{mainremark}.

\smallskip
The inverse matrix of $(g_{ij})_x$ is
\[
(g^{ij})_x = \begin{pmatrix}
-\frac{1}{c^2}+x_2^2 &-x_2 &0\\
-x_2 &1 &0\\
0 &0 &1
\end{pmatrix} \;\;.
\]

We compute the Christoffel symbols (symmetric in the two lower indices)
\[
\varGamma_{ij}^k = \frac{1}{2}\sum_mg^{km}(\partial_ig_{jm}+\partial_jg_{im}-\partial_mg_{ij})
\]
of the metric $g^c_3$ with respect to the coordinates $(x_0,x_1,x_2)$:
\begin{align*}
\varGamma_{00}^0 &= 0  & \varGamma_{00}^1 &= 0  & \varGamma_{00}^2 &= 0\\
\varGamma_{01}^0 &= 0  & \varGamma_{01}^1 &= 0  & \varGamma_{01}^2 &= \tfrac{c^2}{2}\\
\varGamma_{02}^0 &= \tfrac{c^2}{2}x_2  & \varGamma_{02}^1 &= -\tfrac{c^2}{2}  & \varGamma_{02}^2 &= 0\\
\varGamma_{11}^0 &= 0  & \varGamma_{11}^1 &= 0  & \varGamma_{11}^2 &= c^2x_2\\
\varGamma_{12}^0 &= \tfrac{1}{2}(1+c^2x_2^2) & \varGamma_{12}^1 &= -\tfrac{c^2}{2}x_2 & \varGamma_{12}^2 &= 0\\
\varGamma_{22}^0 &= 0  & \varGamma_{22}^1 &= 0  & \varGamma_{22}^2 &= 0
\end{align*}

Now we compute the geodesics. The (maximal) geodesic $\gamma = \sum_{k=0}^2\gamma_k\partial_k$ with $\gamma(0)=p\in\R^3$ and $\gamma'(0)=v\in T_p\R^3=\R^3$ solves for every $k\in\set{0,1,2}$ the equation
\[
\gamma_k''(t) = -\sum_{i,j=0}^2\varGamma_{ij}^k(\gamma(t))\gamma_i'(t)\gamma_j'(t) \;\;.
\]
That is,
\begin{align*}
\gamma_0''(t) &= -2\varGamma_{02}^0(\gamma(t))\gamma_0'(t)\gamma_2'(t) -2\varGamma_{12}^0(\gamma(t))\gamma_1'(t)\gamma_2'(t)\\
&= -c^2\gamma_2(t)\gamma_0'(t)\gamma_2'(t) -(1+c^2\gamma_2(t)^2)\gamma_1'(t)\gamma_2'(t) \;,\\
\gamma_1''(t) &= -2\varGamma_{02}^1(\gamma(t))\gamma_0'(t)\gamma_2'(t) -2\varGamma_{12}^1(\gamma(t))\gamma_1'(t)\gamma_2'(t)\\
&= c^2\gamma_0'(t)\gamma_2'(t) +c^2\gamma_2(t)\gamma_1'(t)\gamma_2'(t) \;,\\
\gamma_2''(t) &= -2\varGamma_{01}^2(\gamma(t))\gamma_0'(t)\gamma_1'(t) -\varGamma_{11}^2(\gamma(t))\gamma_1'(t)^2\\
&= -c^2\gamma_0'(t)\gamma_1'(t) -c^2\gamma_2(t)\gamma_1'(t)^2 \;.
\end{align*}

Let $\omega\define c^2(v_0+p_2v_1)$. If $\omega\neq0$, then the unique solution is
\begin{align*}
\gamma_0(t) &= \frac{v_1v_2}{2\omega^2}\cos(2\omega t) +\frac{v_2^2-v_1^2}{4\omega^2}\sin(2\omega t) -\frac{p_2\omega -v_1}{\omega^2}\big(v_1\sin(\omega t) -v_2\cos(\omega t)\big)\\
&\mspace{20mu}+\Big(\frac{\omega}{c^2} -\frac{v_1^2+v_2^2}{2\omega}\Big)t
+\frac{2p_0\omega^2 -2p_2v_2\omega +v_1v_2}{2\omega^2} \;\;,\\
\gamma_1(t) &= \frac{v_1\sin(\omega t) -v_2\cos(\omega t) +p_1\omega +v_2}{\omega} \;\;,\\
\gamma_2(t) &= \frac{v_2\sin(\omega t) +v_1\cos(\omega t) +p_2\omega -v_1}{\omega} \;\;.
\end{align*}

If $\omega=0$, then the unique solution is
\begin{align*}
\gamma_0(t) &= -\tfrac{1}{2}v_1v_2t^2 -p_2v_1t +p_0 \;\;,\\
\gamma_1(t) &= v_1t +p_1 \;\;,\\
\gamma_2(t) &= v_2t +p_2 \;\;.
\end{align*}

In particular, the metric $g^c_3$ is geodesically complete. The same holds for the product metric $g^c_n$, so the second statement of Remark \ref{mainremark} is true.

\smallskip
No (nondegenerate) $g^c_n$-geodesic with $\omega=0$ is closed, because closedness would obviously imply $v_i=0$ for $i>2$, and $v_1=v_2=0$ and thus also $v_0=0$. Each geodesic with $\omega\neq0$ is the sum of a periodic path and a linear path; hence it is closed if and only if the linear part vanishes, i.e., if and only if $\omega/c^2 -(v_1^2+v_2^2)/(2\omega) = 0$ and $v_i=0$ for all $i>2$, i.e.\ iff $2c^2(v_0+p_2v_1)^2 = v_1^2+v_2^2$ and $v_i=0$ for all $i>2$. The vector $v\in T_p\R^3$ is $g^c_3$-causal iff $-c^2(v_0+p_2v_1)^2 +v_1^2+v_2^2 = -c^2v_0^2 -2c^2p_2v_0v_1 +(1-c^2p_2^2)v_1^2 +v_2^2$ is $\leq0$. So if $v$ were $g^c_3$-causal and tangential to a closed geodesic, then $-c^2(v_0+p_2v_1)^2 +v_1^2+v_2^2 \leq0 = -2c^2(v_0+p_2v_1)^2 +v_1^2+v_2^2$, hence $v_0+p_2v_1=0$ and $v_1^2+v_2^2=0$, i.e.\ $v=0$, a contradiction. Thus there is no closed causal geodesic. But for every $p\in\R^3$, there exist many closed geodesics through $p$: Let $I\subseteq\R$ denote the set consisting of all $v_1\in\R$ with $q(v_1)\define 2c^2(1+p_2v_1)^2 -v_1^2 \geq0$; clearly $I$ is a neighb{\ou}rhood of $0$. For each $v_1\in I$, the geodesic through $p$ in the direction $(1,v_1,v_2)$, where $v_2=\pm\sqrt{q(v_1)}$, is closed. This proves the third statement of Remark \ref{mainremark}.
\end{proof}

\begin{proof}[Proof of the fourth statement in Remark \ref{mainremark} (about existence of closed timelike paths)]{\ }\smallskip\\
Let $w\colon[a,b]\to\R^3$ be a $C^r$ path which is $g^c_3$-timelike. We have to construct for some $\alpha<a$ and $\beta>b$ a $C^1$ extension $w\colon[\alpha,\beta]\to\R^3$ of $w$ which is $g^c_3$-timelike and satisfies $w(\alpha)=w(\beta)$ and $w'(\alpha)=w'(\beta)$. That suffices to prove the claim because every $C^1$ path $w$ which is $C^r$ on some closed interval $I$ can be $C^1$-approximated by $C^r$ paths which are equal to $w$ on $I$, and sufficiently good $C^1$-approximations of timelike paths are timelike.

\smallskip
\emph{First step.} We $C^1$-extend $w$ to an interval $[a_0,b_0]$ with $a_0<a$ and $b_0>b$, such that the extension is still $g^c_3$-timelike and satisfies $w(a_0)=(p_a,0,0)$ and $w(b_0)=(p_b,0,0)$ and $w'(a_0)=w'(b_0)=\nu(\frac{1}{c^2},0,0)$ for some $p_a,p_b\in\R$, $\nu\in\set{1,-1}$.

\smallskip
In order to see that this is possible, we $C^1$-extend the $1$- and $2$-components $w^1,w^2$ of $w$ from $[a,b]$ to an interval $[a_0',b_0']$ such that their values and derivatives at $a_0',b_0'$ are $0$. Now we $C^1$-extend also the $0$-component $w^0$ of $w$ to $[a_0',b_0']$, such that the resulting path $w\in C^1([a_0',b_0'],\R^3)$ is $g^c_3$-timelike. This can in fact be done; we explain the extension to $[a,b_0']$, the extension to $[a_0',b]$ works analogously. For $t\in[a,b_0']$, consider the quadratic polynomial $Q_t\colon\R\to\R$ given by
\[
Q_t(z)\define -c^2z^2 -2c^2w^2(t)(w^1)'(t)z +(1-c^2w^2(t)^2)(w^1)'(t)^2 +(w^2)'(t)^2 \;\;.
\]
We have to choose an extension $v^0\in C^0([a,b_0'],\R)$ of $(w^0)'\in C^0([a,b],\R)$, such that $Q_t(v^0(t))<0$ holds for all $t\in[b,b_0']$. For $t\in[b,b_0']$, we define $F(t)\in\R$ to be the largest $z\in Q_t^{-1}(\set{0})$ if $Q_t^{-1}(\set{0})\neq\leer$, and we define it to be the unique $z$ where $Q_t$ is maximal if $Q_t^{-1}(\set{0})=\leer$. Since $(t,z)\mapsto Q_t(z)$ is continuous, so is $F\colon[b,b_0']\to\R$. We define $f\colon[b,b_0']\to\R$ analogously, just replacing the word ``largest'' by ``smallest''. Since $w'(b)$ is $g^c_3$-timelike, we have $Q_b((w^0)'(b))<0$ and thus $(w^0)'(b)\leq f(b)$ or $F(b)\leq (w^0)'(b)$. If $(w^0)'(b)\leq f(b)$, we define $v^0(t)\define f(t)-f(b)+(w^0)'(b)+b-t$ (which yields $v^0(b)=(w^0)'(b)$ and $v^0(t)<f(t)$ for $t>b$, thus $Q_t(v^0(t))<0$ everywhere). Otherwise we define $v^0(t)\define F(t)-F(b)+(w^0)'(b)-b+t$ (which yields $v^0(b)=(w^0)'(b)$ and $v^0(t)>F(t)$ for $t>b$, thus again $Q_t(v^0(t))<0$ everywhere). We extend $w^0\in C^1([a,b],\R)$ to the interval $[a,b_0']$ by $(w^0)'=v^0$.

\smallskip
We have now obtained a $g^c_3$-timelike path $w\in C^1([a_0',b_0'],\R^3)$ such that $w(a_0)=\!(p_a',0,0)$ and $w(b_0)=(p_b',0,0)$ and $w'(a_0)=(v_a,0,0)$ and $w'(b_0)=(v_b,0,0)$ for some $p_a',p_b'\in\R$, $v_a,v_b\in\R$. Since $g^c_3$ is time-oriented by the vector field $\partial_0$, the fact that the vectors $w'(a_0),w'(b_0)$ are either both future-directed or both past-directed implies $v_av_b>0$. We extend $w$ to an interval $[a_0,b_0]$ by changing $w'$ affinely on $[a_0,a_0']$ and $[b_0',b_0]$ in such a way that $w'$ vanishes nowhere, and that $w'(a_0)=w'(b_0)=\nu(\frac{1}{c^2},0,0)$ for some $\nu\in\set{1,-1}$.

\smallskip
\emph{Second step.} Let $\varphi_c\colon\R^3\to\R^3$ denote the diffeomorphism with $g^c_n = \frac{1}{c^2}\varphi_c^\ast(g^1_n)$ described in Remark \ref{mainremark}, and let $\tilde{w}\define \varphi_c\compose w\in C^1([a,b],\R^3)$. This $\tilde{w}$ is $g^1_3$-timelike (because $w$ is $g^c_3$-timelike) and satisfies $\tilde{w}'(a)=\tilde{w}'(b)=(\nu,0,0)$. In order to prove the fourth statement in Remark \ref{mainremark}, it suffices to extend $\tilde{w}$ to a $g^1_3$-timelike closed path. Since reflections $(x_0,x_1,x_2)\mapsto(q-x_0,x_1,-x_2)$ are $g^1_3$-isometries, we can also assume $\tilde{w}(a)=(0,0,0)$ and $\tilde{w}'(a) =\tilde{w}'(b)=(1,0,0)$ without loss of generality.

\smallskip
Summing up, it suffices to prove for all $p\in\R$ that there exists a $g^1_3$-timelike $C^1$ path $\tilde{w}\colon[\tilde{a},\tilde{b}]\to\R^3$ with $\tilde{w}(\tilde{a})=(p,0,0)$ and $\tilde{w}(\tilde{b})=(0,0,0)$ and $\tilde{w}'(\tilde{a})=\tilde{w}'(\tilde{b})=(1,0,0)$.

\smallskip
\emph{Third step.} Note that $v\in T_x\R^3$ is $g^1_3$-timelike iff $-v_0^2 -2x_2v_0v_1 +(1-x_2^2)v_1^2 +v_2^2 < 0$, i.e.\ iff $(v_0+x_2v_1)^2 > v_1^2+v_2^2$. We choose $T\in\R_{>0}$ so large that $1-0.99(T+1)<\cot(3.14)$. Then the following $C^1$ paths $w_0,\dots,w_{10}$ are all $g^1_3$-timelike, as one can check easily:
\begin{align*}
w_0&\colon[0,1]\to\R^3 &&\text{given by}
&w_0(t) &= \begin{pmatrix}
p+t\\ 0\\ \tfrac{0.99}{2}t^2
\end{pmatrix}\;\;,\\
w_1&\colon[0,T]\to\R^3 &&\text{given by}
&w_1(t) &= \begin{pmatrix}
p+1+t\\ 0\\ 0.99(\tfrac{1}{2}+t)
\end{pmatrix}\;\;,\\
w_2&\colon[0,1]\to\R^3 &&\text{given by}
&w_2(t) &= \begin{pmatrix}
p+T+1+t\\ 0\\ 0.99(\tfrac{1}{2}+T+t-\tfrac{1}{2}t^2)
\end{pmatrix}\;\;,\\
w_3&\colon[0,3.14]\to\R^3 &&\text{given by}
&w_3(t) &= \begin{pmatrix}
p+T+2+\sin(t)\\ 1-\cos(t)\\ 0.99(T+1)
\end{pmatrix}\;\;,\\
w_4&\colon[0,T_1]\to\R^3 &&\text{given by}
&w_4(t) &= \begin{pmatrix}
p+T+2+\sin(3.14)+\cos(3.14)t\\ 1-\cos(3.14)+\sin(3.14)t\\ 0.99(T+1)
\end{pmatrix}\;\;.
\end{align*}
Here we choose $T_1>0$ so large that
\[
-S\define p+2T+6+\pi+2\sin(3.14)-2\cos(3.14)+(\cos(3.14)+\sin(3.14))T_1<0 \;\;.
\]
Let $B^0\define p+T+2+2\sin(3.14)+\cos(3.14)T_1$ and $B^1\define 2-2\cos(3.14)+\sin(3.14)T_1$.

\begin{align*}
w_5&\colon[0,3.14]\to\R^3 &&\text{given by}
&w_5(t) &= \begin{pmatrix}
B^0-\sin(3.14-t)\\ B^1-1+\cos(3.14-t)\\ 0.99(T+1)
\end{pmatrix}\;\;,\\
w_6&\colon[0,1]\to\R^3 &&\text{given by}
&w_6(t) &= \begin{pmatrix}
B^0+t\\ B^1\\ 0.99(T+1-\tfrac{1}{2}t^2)
\end{pmatrix}\;\;,\\
w_7&\colon[0,T]\to\R^3 &&\text{given by}
&w_7(t) &= \begin{pmatrix}
B^0+1+t\\ B^1\\ 0.99(T+\frac{1}{2}-t)
\end{pmatrix}\;\;,\\
w_8&\colon[0,1]\to\R^3 &&\text{given by}
&w_8(t) &= \begin{pmatrix}
B^0+1+T+t\\ B^1\\ 0.99(\frac{1}{2}-t+\frac{1}{2}t^2)
\end{pmatrix}\;\;,\\
w_9&\colon[0,\pi]\to\R^3 &&\text{given by}
&w_9(t) &= \begin{pmatrix}
B^0+2+T+t+\tfrac{1}{2}B^1(1-\cos(t))\\ \tfrac{1}{2}B^1(\cos(t)+1)\\ 0
\end{pmatrix}\;\;,\\
w_{10}&\colon[0,S]\to\R^3 &&\text{given by}
&w_{10}(t) &= \begin{pmatrix}
-S+t\\ 0\\ 0
\end{pmatrix}\;\;.
\end{align*}
Let $b_i>0$ denote the maximum of the domain of $w_i$. Since $w_i(b_i)=w_{i+1}(0)$ and $w_i'(b_i)=w_{i+1}'(0)$ hold for all $i\in\set{0,\dots,9}$, the concatenation of $w_1,\dots,w_{10}$ is a $g^1_3$-timelike path from $(p,0,0)$ to $(0,0,0)$ whose derivatives in the start and end point are $(1,0,0)$.
\end{proof}

\begin{proof}[Proof of Theorem \ref{mainRn} and of the fifth statement in Remark \ref{mainremark}]{\ }

\smallskip
Let $(\partial_0,\dots,\partial_{n-1})$ denote the standard frame of $T\R^n$. We consider the $g^c_n$-orthonormal frame $(e_0,\dots,e_{n-1})$ of $T\R^n$, where $e_0,e_1,e_2$ are given by the formulae in Definition \ref{maindef}, and where $e_i=\partial_i$ for $i\geq3$. The $(n-1)$-plane distribution $H = \spann(e_1,\dots,e_{n-1})$ on $\R^n$ is $g^c_3$-spacelike. It is also twisted, because $[e_1,e_2] = c\,e_0$ vanishes nowhere.

\smallskip
We consider the Riemannian metric $g_V$ on $V\define\spann(\partial_0)$ given by $g_V(\partial_0,\partial_0)=1$, and the Riemannian metric $g_H$ on $H$ which is the restriction of $g^1_n$. Then $(-c^2g_V)\oplus g_H = g^c_n$. By the Key lemma \ref{keylemma}, there exists for every nonempty open set $U\subseteq\R^n$ a $c_U\in\R_{>0}$ such that for all $c\geq k_U$, the set $U$ does not admit a codimension-one foliation none of whose tangent vectors is $g^c_n$-timelike.

\smallskip
Now we compute the Ricci tensor of $g\define g^c_3$. (The result is a special case of the formula in Proposition \ref{Ricformula} below. While we do not spell out the proof there, the computation here is short enough to write down the details.) All Lie brackets $[e_i,e_j]$ except $[e_1,e_2]$ and $[e_2,e_1]$ vanish. For the {\LeviCivita} connection $\nabla$ of $g$, the ``orthonormal Christoffel symbols'' $\Gamma_{ij}^k\define g(\nabla_{e_i}e_j,e_k)$ are given by
\[
2\Gamma_{ij}^k = g([e_i,e_j],e_k) +g([e_k,e_i],e_j) +g([e_k,e_j],e_i) \;\;.
\]
(Note that we used the notation $\varGamma_{ij}^k$ above for the usual Christoffel symbols defined by local coordinates.) Thus all $\Gamma_{ij}^k$ vanish except the following ones:
\begin{align*}
\Gamma_{01}^2 &= \phantom{-}\frac{c}{2} \;\;,  & \Gamma_{10}^2 &= \phantom{-}\frac{c}{2} \;\;,  & \Gamma_{21}^0 &= \phantom{-}\frac{c}{2} \;\;,\\
\Gamma_{02}^1 &= -\frac{c}{2} \;\;,  & \Gamma_{12}^0 &= -\frac{c}{2} \;\;,  & \Gamma_{20}^1 &= -\frac{c}{2} \;\;.
\end{align*}

Since $\nabla_{e_k}e_k=0$ and all functions $\Gamma_{ij}^k$ are constant, we obtain:
\begin{align*}
\Ric_g(e_i,e_j) &= \sum_k\eps_k\Riem_g(e_i,e_k,e_k,e_j)\\
&= \sum_k\eps_k\Big(g(\nabla_{e_i}\nabla_{e_k}e_k,e_j) -g(\nabla_{e_k}\nabla_{e_i}e_k,e_j) -g(\nabla_{[e_i,e_k]}e_k,e_j)\Big)\\
&= -\sum_k\eps_kg(\nabla_{e_k}\nabla_{e_i}e_k,e_j) -\sum_k\eps_kg(\nabla_{[e_i,e_k]}e_k,e_j)\\
&= -\sum_{k,l}\eps_k\eps_lg(\nabla_{e_k}(\Gamma_{ik}^le_l),e_j) -\sum_{k,l}\eps_k\eps_lg([e_i,e_k],e_l)g(\nabla_{e_l}e_k,e_j)\\
&= -\sum_{k,l}\eps_k\eps_l\Gamma_{ik}^l\Gamma_{kl}^j -\sum_{k,l}\eps_k\eps_l\Gamma_{ik}^l\Gamma_{lk}^j +\sum_{k,l}\eps_k\eps_l\Gamma_{ki}^l\Gamma_{lk}^j\\
&= -\sum_{k,l}\eps_k\eps_l\Gamma_{ik}^l\Gamma_{kl}^j -\sum_{k,l}\eps_k\eps_l\Gamma_{il}^k\Gamma_{kl}^j -\sum_{k,l}\eps_k\eps_l\Gamma_{kl}^i\Gamma_{lk}^j\\
&= -\sum_{k,l}\eps_k\eps_l\Gamma_{kl}^i\Gamma_{lk}^j \;\;.
\end{align*}
This yields
\begin{align*}
\Ric(e_0,e_j) &= -\sum_{k,l}\eps_k\eps_l\Gamma_{kl}^0\Gamma_{lk}^j
= -\Gamma_{12}^0\Gamma_{21}^j -\Gamma_{21}^0\Gamma_{12}^j
= -2\delta_{0j}\Gamma_{12}^0\Gamma_{21}^0
= \delta_{0j}c^2/2 \;\;,\\
\Ric(e_1,e_j) &= -\sum_{k,l}\eps_k\eps_l\Gamma_{kl}^1\Gamma_{lk}^j
= \Gamma_{02}^1\Gamma_{20}^j +\Gamma_{20}^1\Gamma_{02}^j
= 2\delta_{1j}\Gamma_{02}^1\Gamma_{20}^1
= \delta_{1j}c^2/2 \;\;,\\
\Ric(e_2,e_j) &= -\sum_{k,l}\eps_k\eps_l\Gamma_{kl}^2\Gamma_{lk}^j
= \Gamma_{01}^2\Gamma_{10}^j +\Gamma_{10}^2\Gamma_{01}^j
= 2\delta_{2j}\Gamma_{01}^2\Gamma_{10}^2
= \delta_{2j}c^2/2 \;\;.
\end{align*}

Hence $\Ric_g$ is with respect to the frame $(e_0,\dots,e_{n-1})$ in each point given by the diagonal matrix $\frac{c^2}{2}\diag(1,1,1,0\dots,0)$. In particular, $\scal_g = c^2/2$. The energy-momentum tensor $T = \Ric -\frac{1}{2}\scal\,g +\Lambda g$ is given by the matrix $\diag(\tfrac{3}{4}c^2-\Lambda,\tfrac{1}{4}c^2+\Lambda,\tfrac{1}{4}c^2+\Lambda,0,\dots,0)$.

\smallskip
For every $v=(v_0,v_1,v_2,\hat{v})\in\R^n$ (where $\hat{v}\in\R^{n-3}$ and $v=\hat{v}+\sum_{i=0}^2v_ie_i$), we have
\[
T(v,v) = \Big(\tfrac{3}{4}c^2-\Lambda\Big)v_0^2 +\Big(\tfrac{1}{4}c^2+\Lambda\Big)(v_1^2+v_2^2) +\Big(-\tfrac{1}{4}c^2+\Lambda\Big)\abs{\hat{v}}^2 \;\;.
\]
The vector $v$ is timelike iff $-v_0^2+v_1^2+v_2^2+\abs{\hat{v}}^2 < 0$. If $\tfrac{3}{4}c^2\geq\Lambda$, we get for each such $v$:
\begin{align*}
T(v,v) &\geq \Big(\tfrac{3}{4}c^2-\Lambda\Big)(v_1^2+v_2^2+\abs{\hat{v}}^2) +\Big(\tfrac{1}{4}c^2+\Lambda\Big)(v_1^2+v_2^2) +\Big(-\tfrac{1}{4}c^2+\Lambda\Big)\abs{\hat{v}}^2\\
&= c^2(v_1^2+v_2^2) +\tfrac{1}{2}c^2\abs{\hat{v}}^2\\
&\geq 0 \;\;.
\end{align*}
Thus $g^c_n$ satisfies the weak energy condition with respect to $\Lambda$ if $\tfrac{3}{4}c^2\geq\Lambda$.

\smallskip
The semi-dominant energy condition holds with respect to $\Lambda$ iff for every timelike vector $v$, the vector $w=-\sharp T(.,v)$ satisfies $g(w,w)\leq0$. In our case $w=(w_0,w_1,w_2,\hat{w})$ with $w_0=-g(e_0,w)=T(e_0,v)=v_0T(e_0,e_0)=(\tfrac{3}{4}c^2-\Lambda)v_0$ and $w_1=-(\tfrac{1}{4}c^2+\Lambda)v_1$ and $w_2=-(\tfrac{1}{4}c^2+\Lambda)v_2$ and $\hat{w}=(\tfrac{1}{4}c^2-\Lambda)\hat{v}$. Thus the semi-dominant energy condition holds iff $v_0^2\geq v_1^2+v_2^2+\abs{\hat{v}}^2$ implies
\[
0 \geq -(\tfrac{3}{4}c^2-\Lambda)^2v_0^2 +(\tfrac{1}{4}c^2+\Lambda)^2(v_1^2+v_2^2) +(\tfrac{1}{4}c^2-\Lambda)^2\abs{\hat{v}}^2 \;\;.
\]
This is true whenever $c^2\geq 4\Lambda$:
\begin{align*}
&(\tfrac{3}{4}c^2-\Lambda)^2v_0^2 -(\tfrac{1}{4}c^2+\Lambda)^2(v_1^2+v_2^2) -(\tfrac{1}{4}c^2-\Lambda)^2\abs{\hat{v}}^2\\
&\geq (\tfrac{3}{4}c^2-\Lambda)^2(v_1^2+v_2^2+\abs{\hat{v}}^2) -(\tfrac{1}{4}c^2+\Lambda)^2(v_1^2+v_2^2) -(\tfrac{1}{4}c^2-\Lambda)^2\abs{\hat{v}}^2\\
&= \Big(\tfrac{9}{16}c^4 -\tfrac{3}{2}\Lambda c^2 -\tfrac{1}{16}c^4 -\tfrac{1}{2}\Lambda c^2\Big)(v_1^2+v_2^2) +\Big(\tfrac{9}{16}c^4 -\tfrac{3}{2}\Lambda c^2 -\tfrac{1}{16}c^4 +\tfrac{1}{2}\Lambda c^2\Big)\abs{\hat{v}}^2\\
&= \Big(\tfrac{1}{2}c^4 -2\Lambda c^2\Big)(v_1^2+v_2^2) +\Big(\tfrac{1}{2}c^4 -\Lambda c^2\Big)\abs{\hat{v}}^2\\
&\geq0 \;\;.
\end{align*}
This proves the fifth statement in Remark \ref{mainremark} and thus, together with the observation at the beginning of the proof, also Theorem \ref{mainRn}.
\end{proof}

\begin{remark}
It is easy to check (and follows from what we prove in later sections) that $g^c_n$ satisfies even the strict dominant energy condition with respect to $\Lambda$ and, moreover, the strict causal convergence condition.
\end{remark}

\begin{remark}
The Weyl tensor of the metric $g^c_3$ vanishes, as for every metric in dimension $3$. Thus the Weyl tensor of $g^c_n$ vanishes, too. In particular, the Petrov type of $g^c_4$ is $0$.
\end{remark}

\section{Global existence of twisted distributions}

In order to prove the theorems \ref{maindom} and \ref{maincob}, we have to show that each manifold $M$ of dimension $n\geq4$ admits an $(n-1)$-plane distribution which is not only twisted on some small open set (as in Lemma \ref{lemma0}) but on a large subset of $M$, preferably on the whole manifold $M$. This can be done via Gromov's h-principle for ample partial differential relations (cf.\ e.g.\ \cite{Spring}, Theorem 4.2; or \cite{GromovPDR}) or, in dimension $\geq5$, via Thom's jet transversality theorem (cf.\ e.g.\ \cite{EliashbergMishachev}, Theorem 2.3.2).

\smallskip
Most of the results we need in later sections have been proved in Chapter 5 of \cite{Nardmann2004}, and we will simply cite them from there. This shortcut is the main reason why we treat in Theorem \ref{maindom} only Lorentzian metrics instead of general pseudo-Riemannian ones: In \cite{Nardmann2004}, where only scalar curvature problems were considered, existence of twisted distributions had to be proved. In the present article, where we deal with Ricci curvature problems, we need distributions whose twistedness satisfies a stronger condition than just being pointwise nonzero; namely, we need (spacelike) distributions $H$ such that for all $x\in M$ and each $w\in T_xM/H_x$, there exist $u,v\in H_x$ with $\Twist_H(u,v)=w$. However, $TM/H$ has rank $1$ in the Lorentzian case, so the latter condition is the same as twistedness then. For pseudo-Riemannian metrics of higher index, we would have to go through the arguments of \cite{Nardmann2004} again, but for the stronger condition instead of twistedness. This would be straightforward but tedious, so let us avoid it here.

\medskip
The definition of the \emph{fine $C^0$-topology} (also known as the \emph{Whitney} or \emph{strong $C^0$-topology}) can be found on p.~9 in \cite{Spring}, for instance. Although the compact-open $C^0$-topology would in principle suffice to prove the theorems in the present article, it is more natural to use the fine $C^0$-topology on the space of distributions in the following sections. (Over a compact manifold, both topologies on the space of distributions are equal.) Its nice property in our context is the following: Every space distribution on a semi-Riemannian manifold $(M,g)$ has a neighb{\ou}rhood with respect to the fine $C^0$-topology all of whose elements are spacelike, too. (On a noncompact manifold, the analogous statement for the compact-open $C^0$-topology is false.)

\smallskip
Chapter 5 of \cite{Nardmann2004} contains the following results which are relevant for the present article:

\begin{theorem} \label{topone}
Let $M$ be a manifold of dimension $\geq5$. Then the set of twisted corank-one distributions on $M$ is dense in $\Distr_{n-1}(M)$ with respect to the fine $C^0$-topology.
\end{theorem}
\begin{proof}
\cite{Nardmann2004}, Theorem 5.3.2.
\end{proof}

\begin{remark}
Although this was not mentioned in \cite{Nardmann2004} and is not relevant for the present article either, we remark that Thom's jet transversality theorem implies that not only $C^0$-denseness but even $C^\infty$-denseness holds in the preceding theorem (because $C^\infty$-generic sections satisfy a certain transversality property and are thus twisted by a simple dimension-counting argument: $n<\frac{1}{2}(n-1)(n-2)$). In the $4$-dimensional case of Theorem \ref{toptwo} below, not only $C^\infty$-denseness but also $C^1$-denseness fails in general, however.
\end{remark}

Recall that the \emph{signature} of a compact oriented $4$-manifold $M$ is the signature of its intersection form on the second homology group $H_2(M;\Z)$; cf.\ e.g.\ Chapter 1 in \cite{GompfStipsicz}. If $M$ admits a line distribution (equivalently: if $M$ admits a Lorentzian metric), then its signature is even; cf.\ \cite{Nardmann2004}, Proposition 5.2.15. When one reverses the orientation of $M$, then the signature changes its sign. Thus the statement ``the signature is divisible by $4$'' does for a connected manifold not depend on the choice of orientation, i.e., it makes sense for compact connected \emph{orientable} $4$-manifolds.

\begin{theorem} \label{toptwo}
Let $M$ be a connected orientable $4$-manifold which is either not compact, or is compact with signature divisible by $4$. Let $H$ be an orientable $3$-plane distribution on $M$, let $\mathscr{U}\subseteq\Distr_3(M)$ be a neighb{\ou}rhood of $H$ with respect to the fine $C^0$-topology. Then $\mathscr{U}$ contains a twisted distribution.
\end{theorem}
\begin{proof}
\cite{Nardmann2004}, Theorem 5.3.3.
\end{proof}

A twisted orientable $3$-plane distribution does not exist on compact connected orientable $4$-manifolds whose signature is not divisible by $4$ (\cite{Nardmann2004}, Propositions 5.2.8 and 5.2.15).

\begin{remark}
The additional assumptions in the $4$-dimensional case of Theorem \ref{maindom} have been made because of the assumptions of the preceding theorem. Their origin is an obstruction-theoretic problem: in the proofs, one has to find a continuous section in a fib{\re} bundle over $M$ whose typical fib{\re} is $S^2$. Proceeding in the usual obstruction-theoretic manner, one triangulates $M$ and tries to construct a section first on the $0$-skeleton, then on the $1$-skeleton, \dots, and finally on the $4$-skeleton, i.e.\ on the whole manifold. Since $S^2$ is $1$-connected, a section exists on the $2$-skeleton. The first obstruction arises in the attempt to extend this section to the $3$-skeleton; the orientability assumptions guarantee that this extension is possible. The signature condition arises in the attempt to extend the section from the $3$-skeleton to the $4$-skeleton.
\end{remark}

It remains to discuss existence of twisted $2$-plane distributions on $3$-manifolds, i.e.\ contact structures (in the general not necessarily cooriented sense) on $3$-manifolds.

\begin{theorem} \label{topthree}
Let $M$ be a $3$-manifold. If $M$ is not orientable, then it does not admit a twisted $2$-plane distribution. If $M$ is orientable, then every connected component of the space of $2$-plane distributions on $M$ contains a twisted one.
\end{theorem}
\begin{proof}
Modulo trivialities (cf.\ \cite{Nardmann2004}, Appendix A.4.1, A.4.3), this follows from Gromov's h-principle theorems in the noncompact case, while the compact case is proved in \cite{Eliashberg1989}.
\end{proof}

\begin{remark}
Although we do not use this fact in the present article, it deserves to be mentioned that in dimensions $3$ and $4$, twisted corank-one distributions are locally essentially unique (because the twisted corank-one distributions are precisely the contact structures in dimension $3$ and precisely the \emph{even-contact structures} in dimension $4$; cf.\ \cite{Nardmann2004}, Appendix A.4): When $M$ is a manifold of dimension $n\in\set{3,4}$ and $H$ is a twisted $(n-1)$-plane distribution on $M$, then every point in $M$ has a neighb{\ou}rhood $U$ on which there exist local coordinates $(x_0,\dots,x_{n-1})$ such that $H$ is on $U$ the kernel of the $1$-form $dx_0+x_2dx_1$. In dimension $3$, this is a special case of Darboux' theorem on contact structures that we mentioned already in Section \ref{C0closed}. In dimension $4$, it is a special case of McDuff's theorem on even-contact structures (Proposition 7.2 in \cite{McDuff1987}).
\end{remark}

\section{The Ricci curvature of stretched metrics}

The following notation will be convenient in our considerations below:

\begin{definition} \label{stretchdef}
Let $V$ be a time distribution on a semi-Riemannian manifold $(M,g)$, let $H$ denote its orthogonal complement. Let $g_V$ denote the Riemannian metric on the vector bundle $V$ which is the restriction of $-g$, let $g_H$ denote the Riemannian metric on $H$ which is the restriction of $g$; thus $g=(-g_V)\oplus g_H$ with respect to the decomposition $TM=V\oplus H$. Then we define $\switch(g,V)$ to be the Riemannian metric $g_V\oplus g_H$ on $M$. For a function $f\in C^\infty(M,\R_{>0})$, we define $\stre(g,f,V)$ to be the semi-Riemannian metric $(-\frac{1}{f^2}g_V)\oplus g_H$ on $M$. (We call the process of replacing $g$ by $\switch(g,V)$ ``\emph{switching $g$ in the direction $V$}'', and we call the process of replacing $g$ by $\stre(g,f,V)$ ``\emph{stretching $g$ in the direction $V$ by the factor $1/f$}''.)

\smallskip
In a context where these data $g,V,f$ are given, we define for each vector $u=u_V+u_H\in V\oplus H = TM$ a vector $\bar{u}\in TM$ by $\bar{u}=fu_V+u_H$. (Note that $\bar{g}\define\stre(g,f,V)$ satisfies $\bar{g}(\bar{u},\bar{u}) = \bar{g}(fu_V,fu_V)+\bar{g}(u_H,u_H) = g(u_V,u_V)+g(u_H,u_H) = g(u,u)$; in particular, $\bar{u}$ is $\bar{g}$-timelike/lightlike/spacelike iff $u$ is $g$-timelike/lightlike/spacelike.)
\end{definition}

Our aim in this section is to describe how the Ricci curvature of $\stre(g,f,V)$ differs from that of $g$. We are particularly interested in the case where the function $f>0$ is a very \emph{small constant}. Lemma \ref{keylemma} implies that $\stre(g,f,V)$ does not admit a space foliation when $H$ is twisted and $f$ is sufficiently small.

\smallskip
The precise formula for the Ricci tensor of $\stre(g,f,V)$ is quite complicated, as it contains many summands. When $f$ is a small constant, however, then one summand dominates all the other ones. This summand involves the twistedness of $H$. The resulting connection between nonexistence of space foliations and Ricci curvature is what we alluded to in the introduction.

\smallskip
Although the dominance of the twistedness term for small constant $f$ is all we have to know for the proofs of the theorems \ref{maindom} and \ref{maincob}, we will write down the complete formula for the Ricci tensor of $\stre(g,f,V)$, in the general semi-Riemannian case and for possibly nonconstant $f$. We are not going to describe here all the geometric objects which occur in the formula. Let me refer you to \cite{Nardmann2004} for basic notation instead; see \ref{geomnotation} for precise references.

\begin{definition}
Let $(M,g)$ be a semi-Riemannian manifold, let $\nabla$ denote the {\LeviCivita} connection of $g$, let $U$ be a distribution on $M$ which has a $g$-orthogonal complement $\bot U$ (e.g.\ a timelike or spacelike distribution). We define a section $\Twist_{g,U}$ in $\Lambda^2(U^\ast)\otimes(\bot U)^\ast$ by
\[
\Twist_{g,U}(u,v,w) \define g(\Twist_U(u,v),w) = g([u,v],w) = g(\nabla_uv,w) -g(\nabla_vu,w) \;\;.
\]
We define a section $\Symst_{g,U}$ in $\Sym^2(U^\ast)\otimes(\bot U)^\ast$ by
\[
\Symst_{g,U}(u,v,w) \define g(\nabla_uv,w) +g(\nabla_vu,w) \;\;.
\]
(The defining expressions of $\Twist_{g,U}$ and $\Symst_{g,U}$ are a priori well-defined for vector fields, and because of their $C^\infty(M,\R)$-linearity a posteriori also for vectors.)
\end{definition}

\begin{notation} \label{geomnotation}
The notations $\divergence_g^U$, \;$\eval{.}{.}_{g,U}$, \;$\laplace_{g,U}^U(f)$ are defined in \cite{Nardmann2004}, \S2.2.1. The tensor field $\Quiem^U_g$ is what has been denoted by $Q^U_g$ in \cite{Nardmann2004}, \S2.2.4.

\smallskip
Note that $\laplace_{g,U}^U(f)$, \;$\Hess_g(f)$, and all terms containing $df$ vanish when the function $f$ is constant.

\smallskip
Let $W_1,\dots,W_k,W$ be vector spaces, let $i,j\in\set{1,\dots,k}$ with $i<j$, let $U$ be a sub vector space of $W_i\cap W_j\cap W$, let $g$ be a symmetric bilinear form on $W$ whose restriction to $U$ is nondegenerate, let $T\in W_1^\ast\otimes\dots\otimes W_k^\ast$. We denote the trace in the $i$th and $j$th index of $T$ over the sub vector space $U\subseteq W_i\cap W_j$ with respect to the metric $g\restrict U$ by $\trace_{g,U}^aT(\dots,a,\dots,a,\dots)$, where $a$ is any free variable which appears here in the $i$th and $j$th index; i.e., using the notation $\eps_\nu\define g(e_\nu,e_\nu)$,
\[
\big(\trace_{g,U}^a T\big)(v_1,\dots,v_{k-2})
= \sum_{\nu=1}^r\eps_\nu T(v_1,\dots,v_{i-1},e_\nu,v_i,\dots,v_{j-2},e_\nu,v_{j-1},\dots,v_{k-2})
\]
for a $g$-orthonormal basis $(e_1,\dots,e_r)$ of $U$. This notation generali{\z}es in an obvious way to the situation where $W_1,\dots,W_k,W$ are vector bundles over a manifold $M$ and $g\in C^\infty(M\ot\Sym^2(W^\ast))$ and $T\in C^\infty(M\ot W_1^\ast\otimes\dots\otimes W_k^\ast)$ are sections.
\end{notation}

\bigskip
\begin{proposition} \label{Ricformula}
Let $V$ be a time distribution on a semi-Riemannian manifold $(M,g)$ of index $q$, let $H$ denote its orthogonal complement, let $f\in C^\infty(M,\R_{>0})$. Then the Ricci tensor of the metric $\bar{g}\define\stre(g,f,V)$ is given by the following formulae (cf.\ \ref{stretchdef} for the notation $\bar{u},\bar{v}$):

\medskip
If $u,v\in H$, then

\[ \begin{split}
&\Ric_{\bar{g}}(\bar{u},\bar{v})\\
&= \Ric_g(u,v)
+\!\tfrac{q}{f}\Hess_g(f)(u,v)
-\tfrac{2q}{f^2}df(u)df(v)
+\tfrac{q(1-f^2)+2f^2}{2f}\eval{\Symst_{g,H}(u,v,.)}{df}_{g,V}\\
&\mspace{20mu}+\tfrac{1}{f}\divergence_g^V(v)df(u)
-\tfrac{f^2(1-f^2)}{4}\trace_{g,V}^{a}\trace_{g,V}^{b}\Twist_{g,V}(a,b,u)\Twist_{g,V}(a,b,v)\\
&\mspace{20mu}+\tfrac{1}{f}\divergence_g^V(u)df(v)
-\tfrac{1-f^2}{2}\trace_{g,V}^{a}\Big(\Quiem_g^H(u,a,v,a)
+\Quiem_g^H(v,a,u,a)\Big)\\
&\mspace{20mu}-\tfrac{1-f^2}{4}\eval{\Symst_{g,H}(u,v,.)}{\divergence_g^H}_{g,V}
-\tfrac{1-f^2}{4}\trace_{g,V}^{a}\trace_{g,V}^{b}\Symst_{g,V}(a,b,u)\Symst_{g,V}(a,b,v)\\
&\mspace{20mu}
-\tfrac{1-f^2}{4}\trace_{g,V}^{a}\trace_{g,V}^{b}\Big(\Symst_{g,V}(a,b,u)\Twist_{g,V}(a,b,v)
+\Symst_{g,V}(a,b,v)\Twist_{g,V}(a,b,u)\Big)\\
&\mspace{20mu}
+\tfrac{1-f^2}{4}\trace_{g,V}^{a}\trace_{g,H}^{b}\Big(\Symst_{g,H}(u,b,a)\Twist_{g,H}(v,b,a)
+\Symst_{g,H}(v,b,a)\Twist_{g,H}(u,b,a)\Big)\\
&\mspace{20mu} -\tfrac{1-f^2}{2f^2}\trace_{g,V}^{a}\trace_{g,H}^{b}\Twist_{g,H}(u,b,a)\Twist_{g,H}(v,b,a) \;\;.
\end{split} \]

\bigskip
If $u,v\in V$, then

\[ \begin{split}
&\Ric_{\bar{g}}(\bar{u},\bar{v})\\
&= f^2\Ric_g(u,v)
+(q-2)f\Hess_g(f)(u,v) +f\laplace_{g,V}^V(f)g(u,v)\\
&\mspace{20mu}+\tfrac{1}{f}\laplace_{g,H}^H(f)g(u,v)
+\tfrac{1}{f}\eval{\divergence_g^V}{df}_{g,H}g(u,v)
+f\eval{\divergence_g^H}{df}_{g,V}g(u,v)\\
&\mspace{20mu}-f\divergence_g^H(v)df(u)
-f\divergence_g^H(u)df(v)
-\tfrac{1}{f}\eval{\Symst_{g,V}(u,v,.)}{df}_{g,H}\\
&\mspace{20mu}-(q-1)\eval{df}{df}_{g,V}g(u,v) -\tfrac{q+1}{f^2}\eval{df}{df}_{g,H}g(u,v)\\
&\mspace{20mu}-\tfrac{(q-2)(1-f^2)}{2f}\eval{\Symst_{g,V}(u,v,.)}{df}_{g,H}
+\tfrac{1-f^2}{4}\eval{\Symst_{g,V}(u,v,.)}{\divergence_g^V}_{g,H}\\
&\mspace{20mu}
-\tfrac{1-f^2}{2}\trace_{g,H}^{a}\Big(\Quiem_g^H(v,a,a,u) +\Quiem_g^H(u,a,a,v)\Big)\\
&\mspace{20mu}-\tfrac{1-f^2}{4}\trace_{g,H}^{a}\trace_{g,V}^{b}\Big(
\Symst_{g,V}(u,b,a)\Twist_{g,V}(v,b,a) +\Symst_{g,V}(v,b,a)\Twist_{g,V}(u,b,a)\Big)\\
&\mspace{20mu}
+\tfrac{1-f^2}{4}\trace_{g,H}^{a}\trace_{g,H}^{b}\Symst_{g,H}(a,b,u)\Symst_{g,H}(a,b,v)\\
&\mspace{20mu}+\tfrac{1-f^2}{4}\trace_{g,H}^{a}\trace_{g,H}^{b}\Big(
\Symst_{g,H}(a,b,u)\Twist_{g,H}(a,b,v) +\Symst_{g,H}(a,b,v)\Twist_{g,H}(a,b,u)\Big)\\
&\mspace{20mu}+\tfrac{f^2(1-f^2)}{2} \trace_{g,H}^{a}\trace_{g,V}^{b}\Twist_{g,V}(u,b,a)\Twist_{g,V}(v,b,a)\\
&\mspace{20mu}
+\tfrac{1-f^2}{4f^2}\trace_{g,H}^{a}\trace_{g,H}^{b}\Twist_{g,H}(a,b,v)\Twist_{g,H}(a,b,u)
\;\;.
\end{split} \]
\newpage

If $u\in V$ and $v\in H$, then
\[ \begin{split}
&\Ric_{\bar{g}}(\bar{u},\bar{v})\\
&= f\Ric_g(u,v)
+\tfrac{1-f^2}{2f}\trace_{g,H}^{a}\Big(\Quiem_g^H(a,a,v,u) -\Quiem_g^H(v,a,a,u)\Big)\\
&\mspace{20mu}+\tfrac{f(1-f^2)}{2}\trace_{g,V}^{a}\!\!\Big(\Quiem_g^V\!(u,a,a,v) -\!\Quiem_g^V\!(a,a,u,v)\!\Big)
-\tfrac{1}{2}\eval{\Symst_{g,H}(v,.,u)}{df}_{g,H}\\
&\mspace{20mu}+(q-1)\Hess_g(f)(u,v)
-\tfrac{q-1}{f}df(u)df(v)
+\divergence_g^V(v)df(u)
-\divergence_g^H(u)df(v)\\
&\mspace{20mu}-\tfrac{(q-1)(1-f^2)+3f^2}{2}\eval{\Twist_{g,V}(u,.,v)}{df}_{g,V}
+\tfrac{(q-1)(1-f^2)+4f^2}{2f^2}\eval{\Twist_{g,H}(v,.,u)}{df}_{g,H}\\
&\mspace{20mu}+\tfrac{1-f^2}{4f}\eval{\Symst_{g,V}(u,.,v)}{\divergence_g^H}_{g,V}
+\tfrac{(1-f^2)(1-2f^2)}{4f}\eval{\Twist_{g,V}(u,.,v)}{\divergence_g^H}_{g,V}\\
&\mspace{20mu}-\tfrac{f(1-f^2)}{4}\eval{\Symst_{g,H}(.,v,u)}{\divergence_g^V}_{g,H}
+\tfrac{(1-f^2)(2-f^2)}{4f}\eval{\Twist_{g,H}(.,v,u)}{\divergence_g^V}_{g,H}\\
&\mspace{20mu}+\tfrac{(1-f^2)^2}{8f}\trace_{g,V}^{a}\trace_{g,H}^{b}\!\! \Big(3\Twist_{g,H}(v,b,a)\Twist_{g,V}(u,a,b) +\!\Symst_{g,H}(v,b,a)\Symst_{g,V}(u,a,b)\Big)\\
&\mspace{20mu}+\tfrac{1-f^4}{8f}\trace_{g,V}^{a}\trace_{g,H}^{b} \Big(\Twist_{g,H}(v,b,a)\Symst_{g,V}(u,a,b) -\Symst_{g,H}(v,b,a)\Twist_{g,V}(u,a,b)\Big) \;\;.
\end{split} \]
\end{proposition}
\begin{proof}
This is a long straightforward computation whose details can be found in \cite{Nardmannformulae}.
\end{proof}

\begin{remark}
In the computation of the preceding formulae, one chooses a $g$-orthonormal frame $(e_1,\dots,e_n)$ of $TM$ such that $(e_1,\dots,e_q)$ is a frame of $V$. Then one can calculate how the ``orthonormal Christoffel symbols'' $\Gamma_{ij}^k\define g(\nabla_{e_i}e_j,e_k)$ change when the metric is stretched. Up to this point, the computation can be found in \cite{Nardmann2004}, \S3.3.1. The result is as follows: Consider the {\LeviCivita} connection $\bar{\nabla}$ of $\bar{g}\define\stre(g,f,V)$ and the $\bar{g}$-orthonormal frame $(\bar{e}_1,\dots,\bar{e}_n)$ which is pointwise defined by \ref{stretchdef}. We write $i:V$ [resp.\ $i:H$] iff $e_i$ is a section in $V$ [resp.\ $H$]. Then $\bar{\Gamma}^k_{ij}\define\bar{g}(\bar{\nabla}_{\bar{e}_i}\bar{e}_j,\bar{e}_k)$ is given by
\[ \begin{split}
\bar{\Gamma}^k_{ij} = \begin{cases}
\Gamma^k_{ij} &\text{if $i,j,k:H$}\\[1ex]
\frac{1}{2}\Big(\frac{1}{f}(\Gamma^k_{ij}-\Gamma^k_{ji}) +f(\Gamma^k_{ij}+\Gamma^k_{ji})\Big) &\text{if $i,j:H$, $k:V$}\\[1ex]
-\frac{1}{2}\Big(\frac{1}{f}(\Gamma^j_{ik}-\Gamma^j_{ki}) +f(\Gamma^j_{ik}+\Gamma^j_{ki})\Big) &\text{if $i,k:H$, $j:V$}\\[1ex]
f\Gamma^k_{ij} +\frac{1}{2}(\frac{1}{f}-f)(\Gamma^i_{kj}-\Gamma^i_{jk}) &\text{if $j,k:H$, $i:V$}\\[1ex]
\Gamma^k_{ij} -\frac{1}{2}(1-f^2)(\Gamma^i_{kj}-\Gamma^i_{jk}) &\text{if $i:H$, $j,k:V$}\\[1ex]
-\frac{1}{2}\Big((\Gamma^j_{ik}+\Gamma^j_{ki}) +f^2(\Gamma^j_{ik}-\Gamma^j_{ki})\Big) -\eps_k\delta_{ik}\frac{1}{f}df(e_j) &\text{if $j:H$, $i,k:V$}\\[1ex]
\frac{1}{2}\Big((\Gamma^k_{ij}+\Gamma^k_{ji}) +f^2(\Gamma^k_{ij}-\Gamma^k_{ji})\Big) +\eps_j\delta_{ij}\frac{1}{f}df(e_k) &\text{if $k:H$, $i,j:V$}\\[1ex]
f\Gamma^k_{ij} -\eps_i\delta_{ik}df(e_j) +\eps_i\delta_{ij}df(e_k) &\text{if $i,j,k:V$}
\end{cases} \;\;.
\end{split} \]
The Riemann tensor of $\bar{g}$ is determined as follows (cf.\ \cite{Nardmann2004}, Formulae 2.2.20):
\[
\Riem_{\bar{g}}(\bar{e}_i,\bar{e}_j,\bar{e}_k,\bar{e}_l) = \partial_{\bar{e}_i}\bar{\Gamma}^l_{jk} -\partial_{\bar{e}_j}\bar{\Gamma}^l_{ik} +\sum_\mu\eps_\mu\Big(\bar{\Gamma}^l_{i\mu}\bar{\Gamma}^\mu_{jk} -\bar{\Gamma}^l_{j\mu}\bar{\Gamma}^\mu_{ik} -(\bar{\Gamma}^\mu_{ij}-\bar{\Gamma}^\mu_{ji})\bar{\Gamma}^l_{\mu k}\Big) .
\]
One can see already at this point that the twistedness of $H$ yields the leading contribution to the curvature of $\bar{g}$ (at the points where is does not vanish) when $f$ is a small constant: The dominant terms are products of terms involving $\frac{1}{f}$, and all corresponding coefficients have the form $\Gamma_{ij}^k-\Gamma_{ji}^k$ with $i,j:H$ and $k:V$, i.e., they are determined by the twistedness of $H$.
\end{remark}

\section{Globali{\z}ation of the dominant energy condition} \label{sectiondominant}

\begin{definition}
Let $(M,g)$ be a semi-Riemannian manifold, let $H$ be a space distribution on $M$, let $V$ denote its $g$-orthogonal complement. We define a section $b_{g,H}$ in the vector bundle $\Sym^2(T^\ast M)\to M$ of symmetric bilinear forms on $TM\to M$ by declaring
\[
b_{g,H}(v,w) \define \begin{cases}
-2\trace_{g,H}^a\trace_{g,V}^b\Twist_{g,H}(v,a,b)\Twist_{g,H}(w,a,b) &\text{if $v,w\in H$}\\
\trace_{g,H}^a\trace_{g,H}^b\Twist_{g,H}(a,b,v)\Twist_{g,H}(a,b,w) &\text{if $v,w\in V$}\\
0 &\text{if $v\in V$ and $w\in H$}
\end{cases} \;\;.
\]
We define another section $\beta_{g,H}$ in $\Sym^2(T^\ast M)\to M$ by $\beta_{g,H}\define b_{g,H}-\frac{1}{2}\trace_g(b_{g,H})g$.

\smallskip
Let $h$ denote the Riemannian metric $\switch(g,V)$ on $M$, let $K$ be a subset of $M$. We say that $(M,g,H)$ is \emph{weak energy nice on $K$} [resp.\ \emph{semi-dominant energy nice on $K$}, resp.\ \emph{causal convergence nice on $K$}] iff there exists a constant $c\in\R_{>0}$ such that for every $x\in K$ and every $g$-causal vector $v\in T_xM$ with $\abs{v}_h\geq1$, we have $\beta_{g,H}(v,v)\geq c$ [resp.\ $-g\big(\sharp(\beta_{g,H}(v,.)), \sharp(\beta_{g,H}(v,.))\big)\geq c$ (where $\sharp\colon T^\ast_xM\to T_xM$ denotes the isomorphism given by $g$), resp.\ $b_{g,H}(v,v)\geq c$]. We say that $(M,g,H)$ is \emph{dominant energy nice on $K$} iff it is weak energy nice on $K$ and semi-dominant energy nice on $K$.
\end{definition}

\begin{proposition} \label{niceprop}
Let $(M,g)$ be a Lorentzian manifold, let $K$ be a compact subset of $M$, let $H$ be a space distribution which is twisted on $K$. Then $(M,g,H)$ is dominant energy nice on $K$ and causal convergence nice on $K$.
\end{proposition}
\begin{proof}
Let $x\in K$. We choose a vector $e_0\in V_x$ with $g(e_0,e_0)=-1$. For $n\define\dim(M)$, let $r\in\N$ denote the number with $2r=n-1$ or $2r+1=n-1$. Since the bilinear form $A\define\Twist_{g,H}(.,.,e_0)\colon H_x\times H_x\to\R$ is skew-symmetric, there exist a vector $\lambda\in\R^r$ and a $g$-orthonormal basis $(e_1,\dots,e_{n-1})$ of $H_x$ such that with respect to this basis, $A$ has the matrix
\[
\begin{pmatrix} && &\lambda_1&&\\ && &&\ddots&\\ && &&&\lambda_r\\
-\lambda_1&& &&&\\ &\ddots& &&&\\ &&-\lambda^r &&& \end{pmatrix}
\mspace{30mu}\text{if $n-1$ is even},
\]
and has this matrix with an additional row and column of zeroes if $n-1$ is odd.

\smallskip
We obtain for $i,j>0$:
\[ \begin{split}
b_{g,H}(e_0,e_0) &= \sum_{k,l=1}^{n-1}A(e_k,e_l)^2 = 2\abs{\lambda}^2 \;\;,\\
b_{g,H}(e_i,e_j) &= 2\sum_{k=1}^{n-1}A(e_i,e_k)A(e_j,e_k)\\
&= \begin{cases}
2\delta_{ij}\lambda_i^2 &\text{if $i,j\leq r$}\\
2\delta_{ij}\lambda_{i-r}^2 &\text{if $r<i,j\leq2r$}\\
0 &\text{else}
\end{cases} \;\;.
\end{split} \]

\smallskip
For a vector $v=\sum_{k=0}^{n-1}v_ke_k\in T_xM$ and $i>0$, we compute
\[ \begin{split}
b_{g,H}(v,v) &= v_0^2b_{g,H}(e_0,e_0) +\sum_{j,k=1}^{n-1}v_jv_kb_{g,H}(e_j,e_k)\\
&= 2v_0^2\abs{\lambda}^2 +2\sum_{k=1}^{r}v_k^2\lambda_k^2 +2\sum_{k=r+1}^{2r}v_k^2\lambda_{k-r}^2 \;\;,\\
b_{g,H}(v,e_0) &= v_0b_{g,H}(e_0,e_0) = 2v_0\abs{\lambda}^2 \;\;,\\
b_{g,H}(v,e_i) &= \sum_{k=1}^{n-1}v_kb_{g,H}(e_k,e_i)
= \begin{cases}
2v_i\lambda_i^2 &\text{if $i\leq r$}\\
2v_i\lambda_{i-r}^2 &\text{if $r<i\leq2r$}\\
0 &\text{if $i=2r+1$}
\end{cases} \;\;.
\end{split} \]

When $v$ is $g$-causal with $\abs{v}_h\geq1$, then $v_0^2-\sum_{i>0}v_i^2\geq0$ and $v_0^2+\sum_{i>0}v_i^2\geq1$, hence $2v_0^2\geq1$. In particular, $b_{g,H}(v,v)\geq\abs{\lambda}^2 = \abs{\Twist_{g,H}}_h^2$. (The equality $\abs{\lambda}^2 = \abs{\Twist_{g,H}}_h^2$ holds because we consider $\Twist_{g,H}$ as a section in $\Lambda^2(H^\ast)\otimes V$. When one considers it as a section in $H^\ast\otimes H^\ast\otimes V$, one gets $\abs{\lambda}^2 = \frac{1}{2}\abs{\Twist_{g,H}}_h^2$.) Because the nonnegative function $\abs{\Twist_{g,H}}_h^2$ vanishes nowhere on the compact set $K$ by assumption, there exists a $c\in\R_{>0}$ such that $b_{g,H}(v,v)\geq c$ holds for all $g$-causal $v\in TM\restrict K$ with $\abs{v}_h\geq1$; i.e., $(M,g,H)$ is causal convergence nice.

\smallskip
Using
\[ \begin{split}
\trace_g(b_{g,H}) &= -b_{g,H}(e_0,e_0)+\sum_{i>0}b_{g,H}(e_i,e_i) = -2\abs{\lambda}^2+4\abs{\lambda}^2 = 2\abs{\lambda}^2 \;\;,
\end{split} \]
we get $\beta_{g,H}(v,v) = b_{g,H}(v,v) -\frac{1}{2}\trace_g(b_{g,H})g(v,v) \geq b_{g,H}(v,v)$ for all $g$-causal $v$. Thus $(M,g,H)$ is weak energy nice.

\smallskip
Since $\sharp(\beta_{g,H}(v,.)) = \sum_{i=0}^{n-1}\eps_i\beta_{g,H}(v,e_i)e_i$, we obtain for $g$-causal $v$:
\[ \begin{split}
&-g\big(\sharp(\beta_{g,H}(v,.)), \sharp(\beta_{g,H}(v,.))\big)\\
&= -\sum_{i,j=0}^{n-1}\eps_i\eps_j\beta_{g,H}(v,e_i)\beta_{g,H}(v,e_j)g(e_i,e_j)
= -\sum_{i=0}^{n-1}\eps_i\beta_{g,H}(v,e_i)^2\\
&= -\sum_{i=0}^{n-1}\eps_i\Big(b_{g,H}(v,e_i)-\abs{\lambda}^2g(v,e_i)\Big)^2\\
&= \Big(3v_0\abs{\lambda}^2\Big)^2 -\sum_{i=1}^{r}\Big(v_i(2\lambda_i^2-\abs{\lambda}^2)\Big)^2 -\sum_{i=r+1}^{2r}\Big(v_i(2\lambda_{i-r}^2-\abs{\lambda}^2)\Big)^2\\
&= 9\abs{\lambda}^4v_0^2
-\sum_{i=1}^{r}\Big(4\lambda_i^4 -4\lambda_i^2\abs{\lambda}^2 +\abs{\lambda}^4\Big)v_i^2
-\sum_{i=r+1}^{2r}\Big(4\lambda_{i-r}^4 -4\lambda_{i-r}^2\abs{\lambda}^2 +\abs{\lambda}^4\Big)v_i^2\\
&\geq 9\abs{\lambda}^4v_0^2 -\sum_{i=1}^{r}\abs{\lambda}^4v_i^2 -\sum_{i=r+1}^{2r}\abs{\lambda}^4v_i^2
\;\geq\; \abs{\lambda}^4\Big(9v_0^2 -\sum_{i=1}^{n-1}v_i^2\Big)
\;\geq\; 8\abs{\lambda}^4v_0^2 \;\;.
\end{split} \]
By a similar argument as above, $(M,g,H)$ is thus semi-dominant energy nice.
\end{proof}

\begin{proposition} \label{niceprop2}
Let $(M,g)$ be a semi-Riemannian manifold, let $H$ be a space distribution on $M$, let $V$ denote its $g$-orthogonal complement, let $K$ be a compact subset of $M$, let $\Lambda\in\R$. If $(M,g,H)$ is
\begin{enumerate} \renewcommand{\labelenumi}{(\alph{enumi})}
\item
weak energy nice
\item
semi-dominant energy nice
\item
dominant energy nice
\item
causal convergence nice
\end{enumerate}
on $K$, respectively, then there exists a constant $\eps_0\in\R_{>0}$ such that for all $\eps\in\ocinterval{0}{\eps_0}$, the metric $\stre(g,\eps,V)$ satisfies the
\begin{enumerate} \renewcommand{\labelenumi}{(\alph{enumi})}
\item
strict weak energy condition with respect to $\Lambda$
\item
strict semi-dominant energy condition with respect to $\Lambda$
\item
strict dominant energy condition with respect to $\Lambda$
\item
strict causal convergence condition
\end{enumerate}
on $K$, respectively.
\end{proposition}
\begin{proof}
Let $h\define\switch(g,V)$. By the third formula from \ref{Ricformula}, there exists a constant $C\in\R_{>0}$ such that all $\eps\in\ocinterval{0}{1}$, all $x\in K$, all $u\in V_x$, and all $v\in H_x$ with $\abs{u}^2_h+\abs{v}^2_h\leq2$ satisfy $\abs{\Ric_{\bar{g}}(\bar{u},\bar{v})} \leq C/\eps$, where $\bar{g}\define\stre(g,\eps,V)$, and where $\bar{u},\bar{v}$ are defined in \ref{stretchdef} (with $f=\eps$).

\smallskip
By the first formula from \ref{Ricformula}, there exists a constant $C_0\in\R_{>0}$ such that all $\eps\in\ocinterval{0}{1}$, $x\in K$, and $u,v\in H_x$ with $\abs{u}^2_h+\abs{v}^2_h\leq2$ satisfy $\abs{\Ric_{\bar{g}}(\bar{u},\bar{v}) -\frac{1}{4\eps^2}b_{g,H}(u,v)} \leq C_0$.

\smallskip
Finally, by the second formula from \ref{Ricformula}, there exists a constant $C_1\in\R_{>0}$ such that all $\eps\in\ocinterval{0}{1}$, $x\in K$, and $u,v\in V_x$ with $\abs{u}^2_h+\abs{v}^2_h\leq2$ satisfy $\abs{\Ric_{\bar{g}}(\bar{u},\bar{v}) -\frac{1}{4\eps^2}b_{g,H}(u,v)} \leq C_1$.

\smallskip
For all $\eps\in\ocinterval{0}{1}$, $x\in K$, $u\in V_x$ and $v\in H_x$ with $\abs{u}^2_h+\abs{v}^2_h\leq1$, the vector $w=u+v$ satisfies therefore
\[ \begin{split}
\Ric_{\bar{g}}(\bar{w},\bar{w}) &= \Ric_{\bar{g}}(\bar{u},\bar{u}) +\Ric_{\bar{g}}(\bar{v},\bar{v}) +2\Ric_{\bar{g}}(\bar{u},\bar{v})\\
&\geq \tfrac{1}{4\eps^2}\big(b_{g,H}(u,u)+b_{g,H}(v,v)\big) -C_0-C_1 -\tfrac{2C}{\eps}\\
&= \tfrac{1}{4\eps^2}b_{g,H}(w,w) -C_0-C_1 -\tfrac{2C}{\eps} \;\;.
\end{split} \]

\smallskip
In order to prove (d), we assume that $(M,g,H)$ is causal convergence nice on $K$; i.e., there exists a constant $c\in\R_{>0}$ such that for all $x\in K$ and every $g$-causal vector $w\in T_xM$ with $\abs{w}_h\geq1$, we have $b_{g,H}(w,w)\geq c$. Whenever $x\in K$, and $\bar{w}\in T_xM$ is $\bar{g}$-causal with $\abs{w}_h=1$, then $\Ric_{\bar{g}}(\bar{w},\bar{w})\geq \tfrac{c}{4\eps^2}-C_0-C_1 -\tfrac{2C}{\eps}$ (because $w$ is $g$-causal). Thus there exists an $\eps_0>0$ such that for all $\eps\in\ocinterval{0}{\eps_0}$ and all $\bar{g}$-causal vectors $\bar{w}$ with $\abs{w}_h=1$, we have $\Ric_{\bar{g}}(\bar{w},\bar{w})>0$. Since every $\bar{g}$-causal vector $\bar{w}'$ has the form $\lambda\bar{w}$, where $\lambda\in\R_{>0}$ and $\bar{w}$ is $\bar{g}$-causal with $\abs{w}_h=1$, we see that $\bar{g}=\stre(g,\eps,V)$ satisfies the strict causal convergence condition when $\eps\leq\eps_0$. This completes the proof of (d).

\smallskip
Let $n\define\dim(M)$. For all $\eps\in\ocinterval{0}{1}$, we have on $K$:
\[ \begin{split}
\abs{\scal_{\bar{g}}-\tfrac{1}{4\eps^2}\trace_gb_{g,H}}
&\leq \sum_{i=0}^{n-1}\abs{\Ric_{\bar{g}}(\bar{e}_i,\bar{e}_i) -\tfrac{1}{4\eps^2}b_{g,H}(e_i,e_i)}
\leq n(C_0+C_1) \;\;.
\end{split} \]
This yields for all $\eps\in\ocinterval{0}{1}$, $x\in K$, and $w,w'\in T_xM$ with $\abs{w}^2_h\leq1$, $\abs{w'}^2_h\leq1$:
\[ \begin{split}
&\big|\Ric_{\bar{g}}(\bar{w},\bar{w}') -\tfrac{1}{2}\scal_{\bar{g}}\bar{g}(\bar{w},\bar{w}') +\Lambda\bar{g}(\bar{w},\bar{w}') -\tfrac{1}{4\eps^2}\beta_{g,H}(w,w')\big|\\
&\leq \big|\Ric_{\bar{g}}(\bar{w},\bar{w}') -\tfrac{1}{4\eps^2}b_{g,H}(w,w')\big|\\
&\mspace{20mu}+\tfrac{1}{2}\big|\scal_{\bar{g}}\bar{g}(\bar{w},\bar{w}') -\tfrac{1}{4\eps^2}(\trace_gb_{g,H})g(w,w')\big|
+\big|\Lambda\bar{g}(\bar{w},\bar{w}')\big|\\
&\leq (C_0+C_1+\tfrac{2C}{\eps}) +\tfrac{n}{2}(C_0+C_1) +\abs{\Lambda} \invdef \tfrac{2C}{\eps}+C_3 \;\;;
\end{split} \]
here we used $\bar{g}(\bar{w},\bar{w}')=g(w,w')$ and $\abs{g(w,w')}\leq \abs{h(w,w')} \leq \abs{w}_h\abs{w'}_h\leq1$.

\smallskip
An argument analogous to the one above in case (d) proves now case (a); i.e., if $(M,g,H)$ is weak energy nice on $K$, then there exists an $\eps_0>0$ such that for all $\eps\in\ocinterval{0}{\eps_0}$, the metric $\bar{g}=\stre(g,\eps,V)$ satisfies the strict weak energy condition with respect to $\Lambda$ on $K$.

\smallskip
Now we prove case (b). Let $(M,g,H)$ be semi-dominant energy nice on $K$, and let $T\define \Ric_{\bar{g}}-\frac{1}{2}\scal_{\bar{g}}\bar{g} +\Lambda\bar{g}$. There exists a constant $C_4\in\R_{\geq0}$ such that $\abs{\beta_{g,H}(\bar{v},\bar{v}')}\leq C_4$ for all $x\in K$ and $v,v'\in T_xM$ with $\abs{v}_h\leq1$, $\abs{v'}_h\leq1$. For all $x\in K$, all $v\in T_xM$ with $\abs{v}_h\leq1$, all $e\in T_xM$ with $\abs{e}_h=1$, and all $\nu\in\set{0,1}$, we obtain:
\[ \begin{split}
\nu T(\bar{v},\bar{e})^2 &= \nu\Big(T(\bar{v},\bar{e}) -\tfrac{1}{4\eps^2}\beta_{g,H}(v,e)\Big)^2 +\tfrac{\nu}{16\eps^4}\beta_{g,H}(v,e)^2\\
&\mspace{20mu}+\tfrac{\nu}{2\eps^2}\beta_{g,H}(v,e)\Big(T(\bar{v},\bar{e}) -\tfrac{1}{4\eps^2}\beta_{g,H}(v,e)\Big)\\
&\leq \tfrac{\nu}{16\eps^4}\beta_{g,H}(v,e)^2 +(\tfrac{2C}{\eps}+C_3)^2 +\tfrac{1}{2\eps^2}C_4(\tfrac{2C}{\eps}+C_3) \;\;;
\end{split} \]
in particular (using a local $g$-orthonormal frame $(e_0,\dots,e_{n-1})$):
\[ \begin{split}
\bar{g}\big(\sharp_{\bar{g}}&(T(\bar{v},.)), \,\sharp_{\bar{g}}(T(\bar{v},.))\big)
= \sum_{i=0}^{n-1}\bar{g}(\bar{e}_i,\bar{e}_i)T(\bar{v},\bar{e}_i)^2\\
&\leq \tfrac{1}{16\eps^4}\sum_{i=0}^{n-1}g(e_i,e_i)\beta_{g,H}(v,e_i)^2 +n(\tfrac{2C}{\eps}+C_3)^2 +\tfrac{n}{2\eps^2}C_4(\tfrac{2C}{\eps}+C_3)\\
&= \tfrac{1}{16\eps^4}g\big(\sharp_g(T(v,.)), \;\sharp_g(T(v,.))\big)  +n(\tfrac{2C}{\eps}+C_3)^2 +\tfrac{n}{2\eps^2}C_4(\tfrac{2C}{\eps}+C_3) \;\;.
\end{split} \]
Now an argument analogous to the one above in case (d) proves that for all sufficiently small $\eps>0$, the metric $\bar{g}=\stre(g,\eps,V)$ satisfies the strict semi-dominant energy condition with respect to $\Lambda$ on $K$. This completes the proof of case (b).

\smallskip
Case (c) of the theorem follows from the cases (a) and (b).
\end{proof}

Now we can prove the following slightly sharpened version of Theorem \ref{maindom}:

\begin{theorem} \label{maindom2}
Let $(M,g)$ be a connected Lorentzian manifold of dimension $n\geq4$, let $K$ be a compact subset of $M$, let $\Lambda\in\R$. If $n=4$, assume that $(M,g)$ is time- and space-orientable, and that either $M$ is noncompact, or compact with intersection form signature divisible by $4$. Let $U$ be a nonempty open subset of $M$. Then there exists a Lorentzian metric $g'$ on $M$ such that
\begin{itemize}
\item
every $g$-causal vector in $TM$ is $g'$-timelike;
\item
$g'$ satisfies the strict causal convergence condition on the set $K$;
\item
$g'$ satisfies the strict dominant energy condition with respect to $\Lambda$ on $K$;
\item
$U$ does not admit any codimension-one foliation none of whose tangent vectors is $g'$-timelike; in particular, $(U,g')$ does not admit a space foliation.
\end{itemize}
\end{theorem}
\begin{proof}
We choose a $g$-space distribution $H'$ on $M$. It has a fine $C^0$-neighb{\ou}rhood in $\Distr_{n-1}(M)$ all of whose elements are spacelike. By Theorem \ref{topone} or Theorem \ref{toptwo}, $M$ admits therefore a twisted $g$-space distribution $H$. Let $V$ denote its $g$-orthogonal complement. For $\eps\in\R_{>0}$, we consider the metric $g_\eps\define\stre(g,\eps,V)$. For all sufficiently small $\eps$, every $g$-causal vector in $TM$ is $g_\eps$-timelike. The Key lemma \ref{keylemma} says that for all sufficiently small $\eps$, the set $U$ does not admit any codimension-one foliation none of whose tangent vectors is $g_\eps$-timelike. Proposition \ref{niceprop} shows that $(M,g,H)$ is dominant energy nice on $K$ and causal convergence nice on $K$. Hence Proposition \ref{niceprop2} implies that for all sufficiently small $\eps$, the metric $g_\eps$ satisfies the strict dominant energy condition with respect to $\Lambda$ on $K$, and it satisfies the strict causal convergence condition on $K$.
\end{proof}

\begin{remarks}
It seems likely that the theorem holds without the additional topological assumptions in the $4$-dimensional case. But one would have to work harder to show this, because Theorem \ref{toptwo} does certainly not hold without these assumptions.

\smallskip
When $M$ is noncompact, then the theorem does not imply the existence of a metric $g'$ which satisfies the two curvature conditions on all of $M$: we cannot choose $K=M$ because of the compactness assumption. This assumption is most likely not necessary either, but getting rid of it would complicate the proof enormously: Proposition \ref{niceprop} could be adapted quite easily to the noncompact case, but that would not help much because the proof of Proposition \ref{niceprop2} fails completely in the noncompact situation. One would have to invent qualitatively new methods in order to deal with that problem.

\smallskip
I have no idea whether one can arrange in Theorem \ref{maindom2} that $g'$ is timelike or lightlike or spacelike geodesically complete (even when we assume that $g$ is geodesically complete).
\end{remarks}

We obtain a weaker statement than Theorem \ref{maindom2} in dimension $3$:

\begin{theorem} \label{maindom3}
Let $(M,g)$ be an orientable Lorentzian $3$-manifold, let $K$ be a compact subset of $M$, let $\Lambda\in\R$, let $U$ be a nonempty open subset of $M$. Then there exists a Lorentzian metric $g'$ on $M$ such that
\begin{itemize}
\item
$g'$ lies in the same connected component of the space of Lorentzian metrics on $M$ as $g$;
\item
$g'$ satisfies the strict causal convergence condition on the set $K$;
\item
$g'$ satisfies the strict dominant energy condition with respect to $\Lambda$ on $K$;
\item
$U$ does not admit any codimension-one foliation none of whose tangent vectors is $g'$-timelike; in particular, $(U,g')$ does not admit a space foliation.
\end{itemize}
If $(M,g)$ admits a spacelike contact structure, then $M$ admits a Lorentzian metric $g'$ which satisfies the properties above and, moreover, has the property that
\begin{itemize}
\item
every $g$-causal vector in $TM$ is $g'$-timelike.
\end{itemize}
\end{theorem}
\begin{proof}
If $(M,g)$ admits a spacelike contact structure $H$, then the proof continues exactly as the proof of Theorem \ref{maindom2}. Otherwise we choose a $g$-space distribution $H'$. Theorem \ref{topthree} shows that the connected component of $\Distr_2(M)$ which contains $H'$ contains also a twisted distribution $H$. We choose any Lorentzian metric $\tilde{g}$ on $M$ which makes $H$ spacelike. By the facts reviewed in \ref{connectedcomponents}, $g$ and $\tilde{g}$ lie in the same connected component of the space of Lorentzian metrics on $M$. Now the proof continues as before.
\end{proof}

\section{Lorentz cobordisms versus weak Lorentz cobordisms} \label{tiplersection}

The term \emph{Lorentz cobordism} was introduced by Yodzis \cite{Yodzis1}, but the discussion of the possibility of ``topology change'' in General Relativity is older; cf.\ e.g.\ the last paragraph in \cite{Reinhart1963}. Geroch observed that topology change can only occur via causality-violating cobordisms (\cite{Geroch1967}, Theorem 2): When there exists a Lorentz cobordism between closed manifolds $S_0$ and $S_1$ which admits no closed timelike curve, then $S_0$ and $S_1$ are diffeomorphic.\footnote{Theorem 2 in \cite{Geroch1967} is false in the form stated there, where it is neither assumed that $S_0$ and $S_1$ are connected, nor that time is inward-directed on $S_0$ and outward-directed on $S_1$. As a counterexample, take any nonempty connected $S$, let $M\define S\times[0,1]\times\set{2,3}$ and $S_0\define S\times\set{0,1}\times\set{2} \cup S\times\set{0}\times\set{3}$ and $S_1\define S\times\set{1}\times\set{3}$. Actually, since Geroch does not assume that $S_0,S_1$ are nonempty, one can even take $M=S\times[0,1]$ and $S_0=S\times\set{0,1}$ and $S_1=\leer$ as a counterexample.

\smallskip The same counterexamples apply to Tipler's comment after Theorem 4 in \cite{Tipler1977}.} Ten years later, Tipler showed that even causality-violating Lorentz cobordisms cannot change the spatial topology provided they satisfy the lightlike convergence condition and some small extra assumption. This extra assumption could be the ``generic condition'', for instance, or the strict lightlike convergence condition. The latter version can be compared nicely to our own results about \emph{weak} Lorentz cobordisms:

\begin{theorem}[Tipler 1977]
Let $n\geq1$, let $S_0,S_1$ be closed $(n-1)$-dimensional manifolds, let $(M,g)$ be a Lorentz cobordism between $S_0$ and $S_1$ which satisfies the strict lightlike convergence condition and has no closed connected component. Then there exists a diffeomorphism $\varphi\colon S_0\times[0,1]\to M$ such that the submanifold $\varphi(\set{x}\times[0,1])$ is $g$-timelike for every $x\in S_0$; in particular, $S_0$ is diffeomorphic to $S_1$.
\end{theorem}
\begin{proof}[Remarks on the proof]
This theorem is not stated explicitly in Tipler's article \cite{Tipler1977}, but it is proved there. Tipler's restriction to $3$-dimensional orientable manifolds in his Theorem 4 is unnecessary. By his Theorem 3, the assumption that $(M,g)$ satisfies the weak energy condition and the generic condition can be replaced by the assumption that the strict lightlike convergence condition holds for $(M,g)$. (Tipler assumes in Theorem 3 the strict weak energy condition [the \emph{ubiquitous energy condition} in his terminology], which implies the strict lightlike convergence condition. The proof shows that only the latter condition is needed.)
\end{proof}

\begin{remark}
In the spirit of Theorem \ref{main0}, one might ask whether for every Lorentzian manifold $(M,g)$ which satisfies a certain Ricci curvature condition on a neighb{\ou}rhood of a closed subset $A\subseteq M$, there exists a Lorentzian metric $g'$ on $M$ which satisfies the condition everywhere and is equal to $g$ on $A$. Tipler's theorem shows that this cannot be true in general (i.e.\ without assumptions on $A$) for the strict lightlike convergence condition, because otherwise one could for instance take nondiffeomorphic $3$-manifolds $S_0,S_1$ and take $A$ to be a small neighb{\ou}rhood of $\mfbd M$ in a Lorentz cobordism $(M,g)$ between $S_0$ and $S_1$ which satisfies the strict lightlike convergence condition on $A$ (such $S_0,S_1,M,g$ do certainly exist).
\end{remark}

In contrast to Tipler's theorem, Theorem \ref{maincob} tells us that the weak Lorentz cobordance relation is hardly restrictive. We state a slight improvement involving the \emph{strict} dominant energy condition:

\begin{theorem} \label{maincob2}
Let $n\geq4$, let $S_0,S_1$ be closed $(n-1)$-dimensional manifolds, let $(M,g)$ be a weak Lorentz cobordism between $S_0$ and $S_1$, let $\Lambda\in\R$. If $n=4$, assume that $M$ is orientable and has no closed connected component. Let $U$ be a nonempty open subset of $M$. Then there exists a weak Lorentz cobordism $(M,g')$ between $S_0$ and $S_1$ such that
\begin{itemize}
\item
every $g$-causal vector in $TM$ is $g'$-timelike;
\item
$(M,g')$ satisfies the strict lightlike convergence condition and the strict dominant energy condition with respect to $\Lambda$;
\item
$U$ does not admit any codimension-one foliation none of whose tangent vectors is $g'$-timelike; in particular, $(U,g')$ does not admit any space foliation.
\end{itemize}
\end{theorem}
\begin{proof}[Proof of Theorem \ref{maincob}]
By shrinking $U$ if necessary, we arrange that $U\cap\mfbd M=\leer$. We extend the Lorentzian manifold-with-boundary $(M,g)$ by an open collar of the boundary, obtaining a manifold $(\bar{M},\bar{g})$. Since every weak Lorentz cobordism is time-orientable by definition, each connected component of $\bar{M}$ satisfies the assumptions of Theorem \ref{maindom2} with $K=M$. Now Theorem \ref{maindom2} yields the result.
\end{proof}

This has the following consequence in dimension $4$:

\begin{corollary}
Let $S_0,S_1$ be oriented closed $3$-manifolds, let $\Lambda\in\R$. Then there exists an oriented weak Lorentz cobordism between $S_0$ and $S_1$ which satisfies the strict causal convergence condition and the strict dominant energy condition with respect to $\Lambda$, and which is an oriented cobordism in the usual sense (i.e., it induces the given boundary orientations).
\end{corollary}
\begin{proof}
Since oriented $3$-manifolds are paralleli{\z}able, Theorem 2 in \cite{Reinhart1963} implies that there exists an oriented cobordism $M$ from $S_0$ to $S_1$ which admits a nowhere vanishing vector field $X$ that is inward-directed on $S_0$ and outward-directed on $S_1$. We remove all closed connected components from $M$ and choose any Lorentzian metric $g$ on $M$ which makes $X$ timelike. Then we apply Theorem \ref{maincob2}.
\end{proof}

Higher dimensions can be discussed similarly, using Reinhart's Theorem 2 and Wall's computation of the oriented cobordism groups in \cite{Wall1960}. It remains to deal with dimension $3$.

\begin{theorem} \label{maincob3}
Let $S_0,S_1$ be closed $2$-manifolds, let $(M,g)$ be an orientable weak Lorentz cobordism between $S_0$ and $S_1$, let $\Lambda\in\R$. Let $U$ be a nonempty open subset of $M$. Then there exists a weak Lorentz cobordism $(M,g')$ between $S_0$ and $S_1$ such that
\begin{itemize}
\item
$g'$ lies in the same connected component of the space of Lorentzian metrics on $M$ as $g$;
\item
$(M,g')$ satisfies the strict lightlike convergence condition and the strict dominant energy condition with respect to $\Lambda$;
\item
$U$ does not admit any codimension-one foliation none of whose tangent vectors is $g'$-timelike; in particular, $(U,g')$ does not admit any space foliation.
\end{itemize}
If $(M,g)$ admits a spacelike contact structure, then $g'$ can be chosen in such a way that
\begin{itemize}
\item
every $g$-causal vector in $TM$ is $g'$-timelike.
\end{itemize}
\end{theorem}
\begin{proof}[Proof of Theorem \ref{maincob}]
By shrinking $U$ if necessary, we arrange that $U\cap\mfbd M=\leer$.  We extend the Lorentzian manifold-with-boundary $(M,g)$ by an open collar of the boundary, obtaining a manifold $(\bar{M},\bar{g})$ which satisfies the assumptions of Theorem \ref{maindom2} with $K=M$. Now Theorem \ref{maindom3} yields the result.
\end{proof}

%%%%%%%%%%%%%%%%%%%%%%%%%%%%%%%%%%%%%%%%%%%%%%%%%%%%%%%%%%%%%%%%%%%%%%%%%%%%%%%%%%%%%%%

%\bibliographystyle{siam}
%\bibliography{bib3copy}

\begin{thebibliography}{10}

\bibitem{Baum}
{\sc H.~Baum}, {\em Spin-{S}trukturen und {D}irac-{O}peratoren \"uber
  pseudoriemannschen {M}annigfaltigkeiten}, vol.~41 of Teubner-Texte zur
  Mathematik [Teubner Texts in Mathematics], BSB B. G. Teubner
  Verlagsgesellschaft, Leipzig, 1981.

\bibitem{BernalSanchez2005}
{\sc A.~N. Bernal and M.~S{\'a}nchez}, {\em Smoothness of time functions and
  the metric splitting of globally hyperbolic spacetimes}, Comm. Math. Phys.,
  257 (2005), pp.~43--50.

\bibitem{Conlon}
{\sc L.~Conlon}, {\em Differentiable manifolds}, Birkh\"auser Advanced Texts:
  Basler Lehrb\"ucher. [Birkh\"auser Advanced Texts: Basel Textbooks],
  Birkh\"auser Boston Inc., Boston, MA, second~ed., 2001.

\bibitem{Eliashberg1989}
{\sc Y.~Eliashberg}, {\em Classification of overtwisted contact structures on
  {$3$}-manifolds}, Invent. Math., 98 (1989), pp.~623--637.

\bibitem{EliashbergMishachev}
{\sc Y.~Eliashberg and N.~Mishachev}, {\em Introduction to the
  {$h$}-principle}, vol.~48 of Graduate Studies in Mathematics, American
  Mathematical Society, Providence, RI, 2002.

\bibitem{Geigessurvey}
{\sc H.~Geiges}, {\em Contact geometry}, in Handbook of differential geometry,
  F.~J.~E. Dillen and L.~C.~A. Verstraelen, eds., vol.~2, 2003, pp.~1--86.

\bibitem{Geroch1967}
{\sc R.~P. Geroch}, {\em Topology in general relativity}, J. Mathematical
  Phys., 8 (1967), pp.~782--786.

\bibitem{GompfStipsicz}
{\sc R.~E. Gompf and A.~I. Stipsicz}, {\em {$4$}-manifolds and {K}irby
  calculus}, vol.~20 of Graduate Studies in Mathematics, American Mathematical
  Society, Providence, RI, 1999.

\bibitem{GromovPDR}
{\sc M.~Gromov}, {\em Partial differential relations}, vol.~9 of Ergebnisse der
  Mathematik und ihrer Grenzgebiete (3) [Results in Mathematics and Related
  Areas (3)], Springer-Verlag, Berlin, 1986.

\bibitem{McDuff1987}
{\sc D.~McDuff}, {\em Applications of convex integration to symplectic and
  contact geometry}, Ann. Inst. Fourier (Grenoble), 37 (1987), pp.~107--133.

\bibitem{Nardmannformulae}
{\sc M.~Nardmann}, {\em Formulae for the curvature of stretched
  semi-{R}iemannian metrics}.
\newblock In preparation.

\bibitem{Nardmann2004}
\leavevmode\vrule height 2pt depth -1.6pt width 23pt, {\em Pseudo-{R}iemannian
  metrics with prescribed scalar curvature}, thesis at the University of
  Leipzig,  (2004), pp.~xviii+158.
\newblock \texttt{arXiv:math.DG/0409435}.

\bibitem{ONeill}
{\sc B.~O'Neill}, {\em Semi-{R}iemannian geometry. {W}ith applications to
  relativity}, vol.~103 of Pure and Applied Mathematics, Academic Press Inc.
  [Harcourt Brace Jovanovich Publishers], New York, 1983.

\bibitem{Reinhart1963}
{\sc B.~L. Reinhart}, {\em Cobordism and the {E}uler number}, Topology, 2
  (1963), pp.~173--177.

\bibitem{Spring}
{\sc D.~Spring}, {\em Convex integration theory. Solutions to the $h$-principle
  in geometry and topology}, vol.~92 of Monographs in Mathematics, Birkh\"auser
  Verlag, Basel, 1998.

\bibitem{Thurston1976}
{\sc W.~P. Thurston}, {\em Existence of codimension-one foliations}, Ann. of
  Math. (2), 104 (1976), pp.~249--268.

\bibitem{Thurston}
\leavevmode\vrule height 2pt depth -1.6pt width 23pt, {\em Three-dimensional
  geometry and topology. {V}ol. 1}, vol.~35 of Princeton Mathematical Series,
  Princeton University Press, Princeton, NJ, 1997.

\bibitem{Tipler1977}
{\sc F.~J. Tipler}, {\em Singularities and causality violation}, Ann. Physics,
  108 (1977), pp.~1--36.

\bibitem{Varela1976}
{\sc F.~Varela}, {\em Formes de {P}faff, classe et perturbations}, Ann. Inst.
  Fourier (Grenoble), 26 (1976), pp.~239--271.

\bibitem{Wall1960}
{\sc C.~T.~C. Wall}, {\em Determination of the cobordism ring}, Ann. of Math.
  (2), 72 (1960), pp.~292--311.

\bibitem{Yodzis1}
{\sc P.~Yodzis}, {\em Lorentz cobordism}, Comm. Math. Phys., 26 (1972),
  pp.~39--52.

\end{thebibliography}

\newcommand{\dummysort}[1]{}

\end{document}